\documentclass[10pt]{article}
\usepackage[a4paper,
centering,scale=0.75]{geometry}
\usepackage{xcolor}
\usepackage{graphicx}
\usepackage{amsmath,amsthm}
\usepackage{amsfonts}
\usepackage{mathrsfs,amsbsy}
\usepackage{mathdots}
\usepackage{fixmath}

\usepackage[shortlabels,inline]{enumitem}

\usepackage{tikz}
\usetikzlibrary{arrows.meta}
\usetikzlibrary{calc}
\usetikzlibrary{datavisualization.formats.functions}

\newtheorem{theorem}{Theorem}[section]
\newtheorem{lemma}{Lemma}[section]

\theoremstyle{definition}
\newtheorem{remark}{Remark}[section]
\newtheorem{example}{Example}[section]
\usepackage{hyperref}
\hypersetup{colorlinks}
\usepackage[capitalize]{cleveref}
\crefname{equation}{}{}
\Crefname{equation}{Equation}{Equations}
\crefrangelabelformat{equation}{(#3#1#4)--(#5#2#6)}

\usepackage{algorithm,algorithmicx,algpseudocode}
\usepackage{subfig}
\usepackage{mathtools}
\DeclarePairedDelimiter\setbasic{\{}{\}}
\DeclarePairedDelimiter\absbasic{|}{|}
\DeclarePairedDelimiter\Nbasic{\|}{\|}
\DeclarePairedDelimiter\innerbasic{\langle}{\rangle}
\newcommand\set[1]{\setbasic*{#1}}
\newcommand\N[1]{\Nbasic*{#1}}
\newcommand\abs[1]{\absbasic*{#1}}
\newcommand\inner[1]{\innerbasic*{#1}}
\renewcommand\Re{\operatorname{Re}}

\def\wtd{\widetilde}
\def\what{\widehat}
\def\ol{\overline}

\DeclareMathOperator{\tr}{tr}

\DeclareMathOperator{\diag}{diag}

\DeclareMathOperator{\HH}{H}
\def\range{{\cal R}}
\def\krylov{{\cal K}}

\DeclareMathOperator{\SD}{SD}
\DeclareMathOperator{\LOCG}{LOCG}

\def\hmz{\hphantom{0.0e-}0}
\DeclareMathOperator{\diff}{d\!}
\DeclareMathOperator{\cheb}{\mathcal{T}}
\DeclareMathOperator{\opL}{\mathscr{L}\!}

\usepackage{multirow}

\allowdisplaybreaks

\title{Local Convergence Behavior of Extended Local Optimal Block Preconditioned Conjugate Gradient Method for Computing Eigenvalues of Hermitian Matrices
}

\author{
	Zhechen Shen\thanks{%
		Department of Mathematical Sciences,
		Tsinghua University,
		Beijing 100084, China.
{\tt shenzc19@mails.tsinghua.edu.cn}.
Supported in part by NSFC-12371380.
	}
	\and
	Xin Liang\thanks{Corresponding author.
		Beijing Institute of Mathematical Sciences and Applications,
		Beijing 101408, China.
{\tt liangxin@bimsa.cn}.
Supported in part by NSFC-12371380.
	}
}
\date{\today}

\begin{document}

\maketitle

\begin{abstract}
This paper provides a comprehensive and detailed analysis of the local convergence behavior of an extended variation of the locally optimal preconditioned conjugate gradient method (LOBPCG) for computing the extreme eigenvalues of a Hermitian matrix.
The convergence rates derived in this work are either obtained for the first time or sharper than those previously established, including those in Ovtchinnikov's work ({\em SIAM J. Numer. Anal.}, 46(5):2567--2592, 2008).
	The study also extends to generalized problems, including Hermitian matrix polynomials that admit an extended form of the Rayleigh quotient.
The new approach used to obtain these rates may also serve as a valuable tool for the convergence analysis of other gradient-type optimization methods.
\end{abstract}

\smallskip
{\bf Keywords.} Extreme eigenvalue, convergence rate, LOBPCG

\smallskip
{\bf AMS subject classifications.} 65F15



\section{Introduction}\label{sec:intro}
The gradient-type optimization methods,
such as the \emph{Steepest Descent} method or the \emph{Conjugate Gradient} method, were originally developed to solve linear systems $Ax=b$ \cite{hestenesS1952methods}, and were later generalized to minimize smooth nonlinear functionals.
For a functional $\varphi(x)$ with the minimal point $x_*$,
The most popular forms of SD and CG read as a two-term recurrence:
\begin{equation*}\label{eq:cg-form}
	x_{i+1}=x_i+\alpha_is_i,\qquad
	s_{i+1}=\nabla \varphi(x_{i+1})+\beta_is_i, i=0,1,\dots, s_0=\nabla \varphi(x_0),
\end{equation*}
where $\nabla \varphi$ is the gradient of $\varphi$.
This recurrence generates a sequence of approximations $x_i$ expected to converge to $x_*$, and a sequence of \emph{conjugate directions} $s_i$.
The stepsize $\alpha_i$ is typically determined by an exact line search on the \emph{search direction} $s_i$ for $\min_\alpha \varphi(x_i+\alpha s_i)$.
Different strategies for selecting $\beta_i$ result in various CG variants, including SD ($\beta_i=0$), Fletcher-Reeves \cite{fletcherR1964function}, Daniel \cite{daniel1967conjugate}, Polak-Ribi\'ere \cite{polakR1969note}, among others.
Another class of variants, introduced by Polyak \cite{polyak1964some} and Takahashi \cite{takahashi1965note}, updates the search direction by historical approximations $x_i$ rather than conjugate directions $s_i$:
\[
	s_{i+1}=\nabla \varphi(x_{i+1})+\beta_ix_i, i=0,1,\dots, s_0=\nabla \varphi(x_0),
\]
This variant is often referred as \emph{Locally Optimal CG}, or momentum acceleration, or heavy ball method.

The extreme (largest or smallest) eigenvalue(s) of a Hermitian matrix admit the variational characterization, or equivalently minimax principles, of the Rayleigh quotient.
Given a Hermitian matrix $A\in \mathbb{C}^{n\times n}$ with eigenvalues $\lambda_1\le\dots\le\lambda_n$, the Rayleigh quotient $\rho(x)$ of a nonzero vector $x\in \mathbb{C}^n$ with respect to $A$ is defined as 
$\rho(x)=\frac{x^{\HH}Ax}{x^{\HH}x}$. 
Here $X^{\HH}$ is the conjugate transpose of a vector or matrix $X$.
The \emph{Fan trace min principle}, one of such minimax principles, states that
\begin{equation*}\label{eq:trace:min}
	\min_{
		X=[x_1,\dots,x_k]\atop
		X^{\HH}X=I_k
	}\tr(X^{\HH}AX)=
	\min_{
		X=[x_1,\dots,x_k]\atop
		X^{\HH}X=I_k
	}\sum_{j=1}^k\rho(x_j)=\sum_{j=1}^k\lambda_j, 	
\end{equation*}
with the minimum attained at some $X$ with $Ax_j=\lambda_jx_j$.
Specially, setting $k=1$ gives 
\begin{equation*}\label{eq:rayleighquotient:min}
	\min_{x\ne 0
	}\rho(x)= \lambda_1. 	
\end{equation*}
Other minimax principles include Courant-Fischer, Wielandt-Lidskii, Amir-Mo\'ez ones.
These minimax principles motivate the use of Rayleigh-Ritz procedures and gradient-type optimization methods
 for computing several extreme eigenvalues and their corresponding eigenvectors, including both conjugate direction variant \cite{bradburyF1966new,perdonG1986exterme,bergamaschiGP1997asymptotic,bergamaschiP2002numerical} and local optimal variant.

In this view, the \emph{locally optimal block preconditioned (extended) conjugate gradient method} (LOBP(e)CG) was developed to solve certain types of eigenvalue problems.
Knyazev \cite{knyazev2001toward} introduced LOBPCG for the generalized Hermitian eigenvalue problem $(A-\lambda B)x=0$, where $A$ is positive definite.
Because of its efficiency, this method has been widely used in various eigenvalue computations.
Unlike the rich convergence results on conjugate direction variants (e.g., global convergence \cite{gilbertN1992global,daiY1999nonlinear}, and convergence rate \cite{crowderW1972linear,cohen1972rate,powell1976some,ritter1980rate}),
up to now, the convergence analysis of LOBPeCG remains incomplete.
To the best of our knowledge, existing results on the estimate for the convergence rate fall into two categories.
For standard Hermitian eigenvalue problems and generalized Hermitian eigenvalue problems $(A-\lambda B)x=0$ with a positive definite $B$, Neymeyr and his co-authors derived the convergence rate of a special form named ``sharp estimate'' for \emph{Preconditioned INverse Vector ITeration} 
and (preconditioned) SD in a series of works \cite{neymeyr2001geometricI,knyazevN2003geometricIII,neymeyrOZ2011convergence,neymeyr2012geometric,argentatiKNOZ2015convergence}, which can be treated as an upper bound of the convergence rate of LOBPCG.
On the other hand,
Ovtchinnikov \cite{ovtchinnikov2008jacobiI,ovtchinnikov2008jacobiII} analyzed the convergence rate of a standard form for LOBPCG applied to standard Hermitian eigenvalue problems and generalized Hermitian eigenvalue problems by 
 establishing its relationship to SD and leveraging the convergence rate of SD by Samokish \cite{samokish1958steepest}.
Along his idea, Benner and Liang \cite{bennerL2022convergence} consider the global convergence and convergence rate of the vector version of LOBPCG on eigenvalue problems associated with definite matrix pencils and hyperbolic quadratic matrix polynomials. 

In this paper, we investigate the local convergence behavior of LOBPeCG by quantitatively comparing the error generated by the approximations $x_{i+1}$ with that generated by the components in $x_{i+1}$ that can be expressed in the previous approximation $x_i$ acted by a polynomial of the Hermitian matrix $A$.
The components can be easily estimated by the polynomial uniformly approximation theory using Chebyshev polynomials.
Using the quantitative relation between the two, we obtain not only new convergence rate for the basic variant of LOBPeCG sharper than the existing ones, but also the completely new ones for complicated variants of LOBPeCG, appearing for the first time.
A brief summary of our contributions is given in \cref{tab:existing-results-and-ours}, after a brief introduction of the algorithm.

The rest of this paper is organized as follows.
In \cref{sec:preliminary}, we present the generic framework of LOBPeCG and introduce two fundamental analyzing tools.
In \cref{sec:a-simple-case,sec:larger-krylov-subspace,sec:more-historical-terms,sec:more-eigenvalues}, we deal with four distinct cases of LOBPeCG:
the basic variant, an extension using larger Krylov subspaces, an extension incorporating more historical approximations, and a version targeting multiple eigenvalues simultaneously.
We conclude in \cref{sec:preconditioners} with a brief discussion of related eigenvalue problems that admit variational characterizations and the role of preconditioners.

\section{Preliminary}\label{sec:preliminary}
Throughout the paper, unless otherwise specified, we focus on the LOBPeCG performing on the Hermitian eigenvalue problem $Ax=\lambda x$ to aim its smallest eigenvalue(s),
where $A=A^{\HH}\in \mathbb{C}^{n\times n}$ is Hermitian,
$\lambda_1\le\dots\le \lambda_n$ are its eigenvalues, associated with eigenvectors $u_1,\dots,u_n$ respectively.
For a rank-$k$ matrix $X\in \mathbb{C}^{n\times k}$,
define the orthogonal projector $P(X)=X(X^{\HH}X)^{-1}X^{\HH}$,
the Rayleigh quotient $\rho(X):=(X^{\HH}X)^{-1/2}X^{\HH}AX(X^{\HH}X)^{-1/2}$,
and
the residual operator $r(X):=AX(X^{\HH}X)^{-1/2}-X(X^{\HH}X)^{-1/2}\rho(X)$.
Clearly $X^{\HH}r(X)=0$.
By $\range(A)$ denote $A$'s column subspace, $\N{A}$ its spectral norm,
while by $\krylov_m(A,X)=\range\left([X,AX,\dots,A^mX]\right)$ denote the length-$(m+1)$ Krylov subspace (of dimension at most $k(m+1)$) generated by $X$ with respect to $A$.

In particular, in the vector case $k=1$, 
for the Rayleigh quotient $\rho(x)=\frac{x^{\HH}Ax}{x^{\HH}x}$,
its gradient $\nabla \rho(x)=\frac{-2}{x^{\HH}x}r(x)$,
and its Hessian $H(x)=(I-2P(x))(A-\rho(x) I)(I-2P(x))$.
Write $F(\lambda)=A-\lambda I$.

Now we present the detailed form $\LOCG(n_b,m_e,m_h)$ of LOBPeCG for computing the $n_b$ smallest eigenvalues in \cref{alg:LOBPeCG}.

\begin{algorithm}[ht]
	\caption{
	$\LOCG(n_b,m_e,m_h)$}\label{alg:LOBPeCG}
	\begin{algorithmic}[1]
		\State Set $X_{-1}=\dots=X_{-m_h}=0$, and choose an initial proper approximation $X_0\in{\mathbb C}^{n\times n_b}$ with $X_0^{\HH}X_0=I_{n_b}$;
		\For{$k=0,1,\dots$}
		\State Construct a nonsingular Hermitian matrix $K_{k}$ as a preconditioner;
		\State Form an orthonormal basis matrix $Z_k$ of the \emph{search subspace}
		\[
			\krylov_{m_e}(K_{k}A,X_k)+\range([X_{k-1},\dots,X_{k-m_h}]);
		\]
		\State Compute the $n_b$ smallest eigenpairs $(\rho_{k+1;j},y_{k;j})$ of the projected problem
		$Z_k^{\HH}F(\lambda)Z_k$, and write diagonal matrix $\rho_{i+1}=\diag(\rho_{i+1;1},\dots,\rho_{i+1;n_b})$, 
		orthonormal matrix $Y_k=[y_{k;1},\dots,y_{k;n_b}]$;
		\State $X_{k+1}=Z_kY_k$.
		\EndFor
	\end{algorithmic}
\end{algorithm}

Here
\begin{itemize}
	\item 
		$n_b\ge 1$ is the number of eigenpairs to compute simultaneously; 
	\item
		$m_h\ge 1$ is the number of historical approximations used for locally optimal purpose;
	\item
		$m_e\ge 1$ indicates the size of the subspace extension so that $m_e+1$ is the length of the Krylov subspace;
	\item
		Ideally $\rho_{k;j}\to \lambda_j$, $X_k\to X_*$ whose columns are the associated eigenvectors respectively.
	\item in Rayleigh-Ritz procedure,  $r(X_{k+1})=AX_{k+1}-X_{k+1}\rho_{k+1},Z_k^{\HH}r(X_{k+1})=0$.
	\item Note that SD is $\LOCG(n_b,1,0)$ and LOBPCG is $\LOCG(n_b,1,1)$.
\end{itemize}

Like all the Rayleigh-Ritz methods for eigenvalues, larger search subspaces yields more accurate approximations.
There are two ways to increase the dimension of the search subspace: one is increasing a higher dimensional Krylov subspace, the other is adding more historical terms (momentums).
Notably, it should be pointed out that adding more historical terms almost brings in no extra calculations.
The reason is given as follows.
The calculations concentrate in generating the search subspace and projecting the matrix into the subspace.
Only considering matrix-vector multiplications of operations $O(n^2)$ other than vector-vector multiplications of operations $O(n)$,
for generating the search subspace, one has to compute $X_k^{(1)}=AX_k, X_k^{(2)}=AX_k^{(2)},\dots,X_k^{(m_e)}=AX_k^{(m_e-1)}$;
for projecting the matrix into the subspace, one has to know $AX_k,AX_k^{(1)},\dots,AX_k^{(m_e)},AX_{k-1},\dots,AX_{k-m_h}$ to obtain the Rayleigh quotient,
in which only $AX_k^{(m_e)}$ has to be computed, for $AX_k=X_k^{(1)},\dots, AX_k^{(m_e-1)}=X_k^{(m_e)}$, and $AX_{k-1},\dots,AX_{k-m_h}$ has been computed in the previous steps.
Hence in each iteration, $(m_e+1)n_b$ matrix-vector multiplications dominate the calculation complexity,
independent of $m_h$.

It is not difficult to show the LOBPeCG's global convergence.
\begin{theorem}\label{thm:global-convergence}
	For the sequences $\set{\rho_{k;j}}$, $\set{x_{k;j}}$ (the $j$th column of $X_k$) generated by $\LOCG(n_b,m_e,m_h)$,
		only one of the following two mutually exclusive situations can occur:
			\begin{enumerate}
				\item 
					for some k, $r(x_{k;j})=0$, and then we have $\rho_{k;j}=\rho_{k+1;j}, x_{k;j}=x_{k+1;j}$;
				\item 
					$\set{\rho_{k;j}}$ is strictly monotonically decreasing, and $\rho_{k;j}\to \what\rho_j$ as $k\to\infty$, while $r(x_{k;j})\neq0$ for all $k$, and no two $x_{k;j}$ are linearly dependent, but $r_{k;j}\to 0$ as $k\to\infty$.
					At this case, for any limit point $\what x_j$ of $\set{x_{k;j}}$, $A\what x_j=\what\rho_j\what x_j$, namely $(\what \rho_j,\what x_j)$ is an eigenpair.
	\end{enumerate}
	Moreover, if $X_*^{\HH}X_0$ is nonsingular, then $\what \rho_j=\lambda_j$ and $\what x_j$ are the associated eigenvectors.
\end{theorem}
\begin{proof}
	With the help of minimax principles, the proof is nearly the same as that of \cite[Theorem 1]{bennerL2022convergence}.
\end{proof}

\Cref{thm:global-convergence} shows that LOBPeCG converges globally, but provides no information on its convergence rate.
In the following, we try to obtain such a rate, starting from the simplest case $\LOCG(1,1,1)$ and gradually extending to the most general case $\LOCG(n_b,m_e,m_h)$. 
A list of existing estimate of the rates and our obtained results is given in \cref{tab:existing-results-and-ours}.

\begin{table}[ht]
	\centering
	\begin{tabular}{|c|cl|c|rc|}
		\hline
		$\LOCG$ & \multicolumn{2}{|c|}{existing results}  && \multicolumn{2}{|c|}{our work}  \\
		\hline
		$(1,m_e,0)$   & \cite{samokish1958steepest,golubY2002inverse}& $C_{m_e}$ & & ---& \\
		$(n_b,m_e,0)$ & \cite{neymeyrOZ2011convergence}& $*$ & &  --- &\\
		$(1,1,1)$   & \cite{ovtchinnikov2008jacobiI}& $\frac{C_{1}}{2-C_{1}}$ &$>$ &$\sqrt{\chi(1,C_{1})}$ &\cref{thm:convergence-rate:111} \\
		$(n_b,1,1)$ & \cite{ovtchinnikov2008jacobiII}& $\frac{C_{1}}{2-C_{1}}$ &$>$& $\sqrt{\chi(\sigma,C_{1})}$ &\cref{thm:convergence-rate:m11} \\
		$(1,m_e,1)$   & \cite{bennerL2022convergence}& $\frac{C_{m_e}}{2-C_{m_e}}$&$>$ & $\sqrt{\chi(1,C_{m_e})}$ &\cref{thm:convergence-rate:1m1} \\
		$(1,1,m_h)$   & &--- && $\sqrt[m_h+1]{\chi(\sigma,C_1)}$ &\cref{thm:convergence-rate:11m} \\
		$(1,m_e,m_h)$ & &---& & $\sqrt[m_h+1]{\chi(\sigma,C_{m_e})}$ &\cref{thm:convergence-rate:1mm} \\
		\hline
	\end{tabular}
	\scriptsize
	\begin{itemize*}
		\item $C_{m_e}=\left(\cheb_{m_e}(\frac{{\kappa}+1}{{\kappa}-1})\right)^{-2},
	\kappa=\frac{\lambda_n-\lambda_1}{\lambda_2-\lambda_1}$,
	$\cheb_{m_e}(\cdot)$ is the $m_e$-degree Chebyshev polynomial of the first kind;
	\end{itemize*}

	\begin{itemize*}
			\item $*$ is complicated;
	\item $>$ implies the inequality;
	\item $\chi(\cdot,\cdot)$ is defined in \cref{eq:chi-psi};
		\item $\sigma$ varies in different terms;
	\end{itemize*}

	\begin{itemize*}
	\item high-order terms are omitted.
	\end{itemize*}
	\caption{Existing results and our work on the estimate of the asymptotic convergence rates}
	\label{tab:existing-results-and-ours}
\end{table}

The key to the results is the relations between the quantities of two successive iterations, described in the following lemma, which appeared as \cite[Theorem~2.2]{bennerL2022convergence}.
\begin{lemma}[{\cite[Theorem~2.2]{bennerL2022convergence}}]\label{cor:-cite-theorem-2-2-bennerl2022convergence-}
	Let vectors $x\ne0,r(x)\ne 0,v\ne0$, and $W\in \mathbb{C}^{n\times m}$, which satisfy $v^{\HH}r(x)\ne0, W^{\HH}r(x)=0$.
	Suppose that $[x, v,W] $ has full column rank,
	and
	$( \alpha_+, b_+)$ is a stationary point of the function $\tr\rho\left(x+ (I-P(x))[\alpha v+Wb] \right)$.
	Write
	\[
		d =   (I-P(x))(\alpha_+v +Wb_+),\quad
		x_+ = x + d.
	\]
	Then, for the nontrivial case that $\alpha_+\ne0$,  
	\begin{subequations}
	\begin{align}
		\label{eq:thm:linesearch:r_opt-perp:1}
		r(x_+)&\perp\range([x,v,W,d]),
		\\
		\alpha_+ 
		&= -\frac{d^{\HH}F(\rho(x_+))d}{r(x)^{\HH}v}, \label{eq:thm:linesearch:alpha:1}
		\\ \rho(x_+)-\rho(x) &= \frac{r(x)^{\HH}v}{x^{\HH}x}\alpha_+ =-\frac{d^{\HH}F(\rho(x_+))d}{x^{\HH}x},
		\label{eq:thm:linesearch:rho:1}
		\\ r(x_+)-r(x) &=\check F(x)d,  \label{eq:thm:linesearch:r:1}
		\\ b_+&= -\alpha_+[W^{\HH}\check F(x)W]^{-1}W^{\HH}\check F(x)v,
		\label{eq:thm:linesearch:beta:1}
	\end{align}
	\end{subequations}
	where $\check F(x)=(I-P(x))F(\rho(x_+))(I-P(x))$,
provided $W^{\HH}\check F(x)W$ is nonsingular.
\end{lemma}
\Cref{cor:-cite-theorem-2-2-bennerl2022convergence-} holds for the particular case $n_b=1$.
We also established an analogue for the general case $n_b\ge 1$ in \cref{thm:linesearch}, provided in \cref{sec:appendix}, due to its technical complexity and the fact that it is not directly invoked in the proofs of the main results below.
\section{Main results}\label{sec:main-results}
	We now present our main results. In order to promptly establish the framework that we used to obtain the asymptotic convergence rates,
	here we only present all the rates that we obtained and a concise proof for the simplest case.
	The details of this concise proof, namely those for \cref{eq:beta,eq:comparison,eq:get-rid-of-phi}, and the complete proofs of the other cases are deferred to \cref{sec:proofs}.

\subsection{The simplest case $\LOCG(1,1,1)$}\label{sec:a-simple-case}
The simplest case of genuine LOBPeCG methods is $\LOCG(1,1,1)$.
For this case, the search subspace is $\range([x_k,r_k,x_{k-1}])$, and the goal is to compute the smallest eigenvalue $\lambda_1$.

At the $k$th iteration, write $r_k=r(x_k),P_k=P(x_k)$, and then the orthogonal basis of the search subspace is $x_k, r_k, (I-P_k)x_{k-1}$ for $x_k\perp r_k, r_k\perp x_{k-1}$.
Since the Rayleigh quotient is invariant under scaling, only the direction of $x_k,x_{k-1},x_{k+1}$ are relevant.
Therefore, we may assume $\N{x_k}=1$ and rewrite the iteration as 
\[
	x_{k+1}=x_k+\alpha_kr_k+\beta_k(I-P_k)x_{k-1},
\]
where $(\alpha_k,\beta_k)$ minimizes the function $\varphi(\alpha,\beta)=\rho(x_k+\alpha r_k+\beta(I-P_k)x_{k-1})$.

Write $\rho_k=\rho(x_k),\varepsilon_k=\rho_k-\lambda_1,d=x_{k+1}-x_k,F_k=A-\rho_k I,\check{F}_{k+1}=(I-P_k)^{\HH}F_{k+1}(I-P_k)$,
$\delta_k=\varepsilon_k-\varepsilon_{k+1}=\rho_k-\rho_{k+1}$.
Then $x_k^{\HH}d_k=0$.
Write $\varphi_k=\angle(x_k,x_{k+1})$, and then
$\N{d_k}^2=\tan^2\varphi_k,\N{d_{k-1}}^2=\sin^2\varphi_{k-1}$,
$\N{x_{k+1}}^2=\sec^2\varphi_k, \N{x_{k-1}}^2=\cos^2\varphi_{k-1}$.

For ease, we drop the subscript $\cdot_k$ and replace $\cdot_{k\pm 1}$ by $\cdot_{\pm}$.
With the help of \cref{cor:-cite-theorem-2-2-bennerl2022convergence-},
the calculation in \cref{sec:proofs} tells us
\begin{equation}\label{eq:beta}
	\alpha=\frac{\delta }{-r ^{\HH}r },\qquad
	\beta =-\frac{\delta }{\delta_-+\delta \sin^2\varphi_-}
	.
\end{equation}

As is well known (see, e.g., \cite{samokish1958steepest,golubY2002inverse}), the SD method $\LOCG(1,1,0)$ enjoys
\begin{equation}\label{eq:rate:SD}
	\varepsilon_+\le C^{\SD}\varepsilon ,
	\quad C^{\SD}= 
	{\cheb_1(\Delta^{-1})}^{-2}+O(\sqrt{\varepsilon })<1,
	\quad
	\Delta=\frac{{\kappa}-1}{{\kappa}+1},
	\;
	\kappa=\frac{\lambda_n-\lambda_1}{\lambda_2-\lambda_1}
	,
\end{equation}
where $\cheb_1(t)=t$ is the first-degree Chebyshev polynomial of the first kind, and $\kappa$  is the condition number of Hermitian eigenvalue problems.

Now we try to compare the convergence behaviors of $\LOCG(1,1,1)$ and SD method $\LOCG(1,1,0)$.
Write $\inner{z,y}_*:=z^{\HH}(A-\lambda_1I)y$, which induces a seminorm $\N{y}_*^2:=\inner{y, y}_*$.
Clearly $\N{x }_*^2=\varepsilon $,
so we expect a relation between $\N{x_+}_*^2
$ and $\N{x +\alpha (A-\lambda_1I)x }_*^2$, 
which is obtained after the calculations in \cref{sec:proofs}:
\begin{equation}
	\N{x +\alpha (A-\lambda_1I)x }_*^2
	=(\alpha \varepsilon_++\alpha ^2\varepsilon ^2)\varepsilon +\alpha \varepsilon_+^2
	+\varepsilon_+
	+\frac{\delta ^2}{\delta_-+(\delta_-+\delta )\tan^2\varphi_-}
	.
	\label{eq:comparison}
\end{equation}

Suppose $\rho \le \rho_-\le \rho_0< \lambda_2$.
First $0\le$ $\varepsilon -\delta =\varepsilon_+\le \varepsilon_+^{\SD}\le C^{\SD} \varepsilon 
\Rightarrow
\varepsilon\ge\delta \ge (1-C^{\SD})\varepsilon .$
Then by the calculations in \cref{sec:proofs}, we have
\begin{equation}\label{eq:get-rid-of-phi}
	\frac{\delta ^2}{\delta_-+(\delta_-+\delta )\tan^2\varphi_-}=
	\frac{\delta ^2}{\delta_-}+O(\delta^2).
\end{equation}
which gives $\N{x +\alpha (A-\lambda_1I)x }_*^2
= \varepsilon_+ +\frac{\delta ^2}{\delta_-}
+O(\varepsilon ^2).$
Let $\N{x +\alpha (A-\lambda_1I)x }_*^2= C\N{x }_*^2=C\varepsilon $, 
and write $\eta =\frac{\varepsilon_+}{\varepsilon }$.
we have
\begin{equation}\label{eq:basic-relation}
	\eta +\frac{\eta_-(1-\eta )^2}{1-\eta_-}= C+O(\varepsilon )=:\wtd C.
\end{equation}

	We need to build an upper bound of $\eta$ and $\eta\eta_-$ from \cref{eq:basic-relation} to reveal the $1$-step and $2$-step asymptotic convergence rate.
	From \cref{eq:basic-relation}, it is clear that $\eta\le \wtd C, \frac{\eta_-(1-\wtd C)^2}{1-\eta_-}\le\wtd C$ and then $\eta_-\le\frac{\wtd C}{1-\wtd C+\wtd C^2}<1$ for $0<\wtd C<1$.
On one hand,
solving the two-variable optimization problem
\[
	\max\; xy, \qquad\text{s.t.}\;
	\begin{cases}
		0\le x\le \wtd C,\quad 0\le y\le \frac{\wtd C}{1-\wtd C+\wtd C^2},& \\
		x+\dfrac{y(1-x)^2}{1-y}=\wtd C,&
	\end{cases}
\]
produces $\eta \eta_-\le \frac{2+\wtd C-2\sqrt{1+\wtd C}}{3-2\sqrt{1+\wtd C}}$, in which the equality is attained at $\eta = \frac{1+\wtd C-(1-\wtd C)\sqrt{1+\wtd C}}{3-\wtd C}$.
Write
\begin{equation}\label{eq:chi-psi}
	\chi(\sigma,\wtd C)= \frac{2\sigma+\wtd C-2\sqrt{\sigma(\sigma+\wtd C)}}{2\sigma+1-2\sqrt{\sigma(\sigma+\wtd C)}},
	\quad
\omega(\sigma,\wtd C)= \frac{\sigma+\wtd C-(1-\wtd C)\sqrt{\sigma(\sigma+\wtd C)}}{1+\sigma(2-\wtd C)},
\end{equation}
and then $\eta\eta_-\le\chi(1,\wtd C)$ attained at $\eta=\omega(1,\wtd C)$.
On the other hand, writing a degree-$1$ polynomial $p_1(\lambda)=1+\alpha \lambda$ and $A=\sum\limits_{i=1}^n \lambda_iu_iu_i^{\HH}$,
\begin{align*}
	C
	&=\frac{\N{x +\alpha (A-\lambda_1I)x }_*^2}{\N{x }_*^2}
	=\frac{\N{p_1(A-\lambda_1I)x }_*^2}{\N{x }_*^2}
	=\frac{\sum_{i\ge 2}(\lambda_i-\lambda_1)\abs{p_1(\lambda_i-\lambda_1)\xi_i}^2}{\sum_{i\ge 2}(\lambda_i-\lambda_1)\abs{\xi_i}^2}
	\\&\le \min_{p_1(0)=1,\deg p_1\le 1}\max_{\lambda\in [\lambda_2-\lambda_1,\lambda_n-\lambda_1]}\abs{p_1(\lambda_i-\lambda_1)}^2
	= 
	{\cheb_1(\Delta^{-1})}^{-2},
\end{align*}
where $\cheb_1(t)=t$ is the first-degree Chebyshev polynomial of the first kind.

We summarize the findings above in \cref{thm:convergence-rate:111}.
\begin{theorem}\label{thm:convergence-rate:111}
	Suppose $\lambda_1\le\rho_0<\lambda_2$.
	For the sequences $\set{\rho_{k}}$, $\set{x_{k}}$  produced by $\LOCG(1,1,1)$,
	\begin{subequations}
	\begin{align}\label{eq: thm:convergence-rate:111}
			\rho_{k+1}-\lambda_1&\le \chi(1,\wtd C)(\rho_{k-1}-\lambda_1), \quad \text{and}
		\\\label{eq: thm:convergence-rate:111:omega} \rho_{k+1}-\lambda_1&\le \omega(1,\wtd C)(\rho_k-\lambda_1) \quad\text{if the equality in \cref{eq: thm:convergence-rate:111} holds},
	\end{align}
	\end{subequations}
	where $\chi(\cdot,\cdot),\omega(\cdot,\cdot)$ are as in \cref{eq:chi-psi} and $\wtd C\le \Delta^2+O(\rho_k-\lambda_1)$.
\end{theorem}

	Note that for $0< C<1$,
	\begin{align*}
		\chi(1, C) &= \frac{C^2}{1+{(1-C)}\left(1+2\sqrt{1+C}\right)}
		< \frac{C^2}{1+(1-C)(3-C)}
		= \left(\frac{C}{2-C}\right)^2
		< C^2
		,
		\\
		\omega(1,C)&=\frac{C}{1+(1-C)(1+C)^{-1/2}}<C
		.
	\end{align*}
Comparing the convergence rate of $\LOCG(1,1,1)$ in \cref{thm:convergence-rate:111} with the convergence rate of $\LOCG(1,1,0)$ in \cref{eq:rate:SD}, we see $\wtd C\le C^{\SD}+O(\rho_k-\lambda_1)$.
Since for $0<\wtd C<1$,
\[
	\text{$2$-step asymptotic convergence rate of $\LOCG(1,1,1)$}\; = \chi(1,\wtd C) < \wtd C^2  \approx \;\text{that of SD},
\]
$\LOCG(1,1,1)$ converges faster than SD method $\LOCG(1,1,0)$ as we expect.
On the other hand, even if the $2$-step convergence rate attains the found upper bound, the $1$-step convergence rate of $\LOCG(1,1,1)$, namely $\omega(1,\wtd C)$, is still better than that of SD method $\LOCG(1,1,0)$, namely $\wtd C$.
Moreover, omitting high-order terms, $(\frac{C^{\SD}}{2-C^{\SD}})^2$ is also the upper bound of the $2$-step convergence rate of $\LOCG(1,1,1)$ estimated by Ovtchinnikov \cite{ovtchinnikov2008jacobiI} (or see \cite{bennerL2022convergence}).
Since for $0<\wtd C<1$,
\[
	\text{our $2$-step asymptotic convergence rate}\; = \chi(1,\wtd C) < (\tfrac{\wtd C}{2-\wtd C})^2  \approx \;\text{Ovtchinnikov's},
\]
here we obtain a sharper estimate of the convergence rate for $\LOCG(1,1,1)$.

\subsection{Larger Krylov subspace: $\LOCG(1,m_e,1)$}\label{sec:larger-krylov-subspace}
For this case, the search subspace is $\range([x_k, r_k, F_kr_k, \dots, F_k^{m_e-1}r_k, x_{k-1}])$.
Since $x_k\perp r_k, r_k\perp x_{k-1}, r_k\perp F_{k-1}^tr_{k-1}$ for $t=0,1,\dots,m_e-1$, 
we have $r_k^{\HH}F_k^tr_{k-1}=r_k^{\HH}(F_{k-1}+\delta_{k-1}I)^tr_{k-1}=0, 
r_k^{\HH}F_k^jx_{k-1}=r_k^{\HH}(F_{k-1}+\delta_{k-1}I)^jx_{k-1}=0.$
Thus,
the orthogonal basis of the search subspace is $x_k, r_k, (I-P_k)x_{k-1}, (I-P_k)(I-Q_k)W_k$ where $Q_k=P(r_k)$, $W_k$ is an orthogonal basis of $\range([F_kr_k, \cdots, F_k^{m_e-1}r_k])$.
So $r_{k-1}^{\HH}W_k=x_{k-1}^{\HH}W_k=0$.
The $k$th iteration is rewritten as
\[
	x_{k+1}=x_k+\alpha_k[r_k+(I-P_k)(I-Q_k)W_kc_k]+\beta_k(I-P_k)x_{k-1},
\]
where $(\alpha_k,\beta_k,c_k)$ minimizes the function $\varphi(\alpha,\beta,c)=\rho(x_k+\alpha q_k+\beta(I-P_k)x_{k-1})$, where 
\begin{equation}\label{eq:qk}
	q_k =r_k+(I-P_k)(I-Q_k)W_kc.
\end{equation}
(Once $r_k\ne 0$, the function has to decrease along the negative gradient direction, which implies $\alpha_k\ne 0$.)

By following the same procedure used to derive the asymptotic convergence rate of $\LOCG(1,1,1)$, we obtain that of $\LOCG(1,m_e,1)$.
\begin{theorem}\label{thm:convergence-rate:1m1}
	Suppose $\lambda_1\le\rho_0<\lambda_2$.
	For the sequences $\set{\rho_{k}}$, $\set{x_{k}}$  produced by $\LOCG(1,m_e,1)$,\!\!
	\begin{subequations}
	\begin{align}\label{eq: thm:convergence-rate:1m1}
			\rho_{k+1}-\lambda_1&\le \chi(1,\wtd C)(\rho_{k-1}-\lambda_1),  \quad \text{and}
			\\\label{eq: thm:convergence-rate:1m1:omega} \rho_{k+1}-\lambda_1&\le \omega(1,\wtd C)(\rho_k-\lambda_1) \quad\text{if the equality in \cref{eq: thm:convergence-rate:1m1} holds},
	\end{align}
	\end{subequations}
	where $\chi(\cdot,\cdot),\omega(\cdot,\cdot)$ are as in \cref{eq:chi-psi} and $\wtd C\le \left(\cheb_{m_e}(\Delta^{-1})\right)^{-2}+O(\rho_k-\lambda_1)$.
	Here $\cheb_{m_e}(t)$ is the $m_e$-th degree Chebyshev polynomial of the first kind.
\end{theorem}
This theorem demonstrates that $\LOCG(1, m_e, 1)$ improves upon the convergence of the eSD method $\LOCG(1, m_e, 0)$. 

\subsection{More historical terms: $\LOCG(1,1,m_h)$}\label{sec:more-historical-terms}
For this case, the search subspace is $\range([x_k,r_k,x_{k-1},\cdots,x_{k-{m_h}}])$.
Let $W_k=[x_{k-1},\dots,x_{k-m_h}]$, and then $r_{k+1}^{\HH}W_k=r_k^{\HH}W_k=0$.
The $k$th iteration is rewritten as
\[
	x_{k+1}=x_k+\alpha_kr_k+(I-P_k)W_kb_k,
\]
where $(\alpha_k,b_k)$ minimizes the function $\varphi(\alpha,b)=\rho(x_k+\alpha r_k+(I-P_k)W_kb)$.

By following the same procedure used to derive the asymptotic convergence rate of $\LOCG(1,1,1)$, we obtain that of $\LOCG(1,1,m_h)$.
\begin{theorem}\label{thm:convergence-rate:11m}
	Suppose $\lambda_1\le\rho_0<\lambda_2$.
	For the sequences $\set{\rho_{k}}$, $\set{x_{k}}$  produced by $\LOCG(1,1,m_h)$,\!\!
	\begin{subequations}
	\begin{align}\label{eq: thm:convergence-rate:11m}
		\rho_{k+1}-\lambda_1&\le \chi(\sigma,\wtd C) 
			(\rho_{k-m_h}-\lambda_1), \quad \text{and}
			\\\label{eq: thm:convergence-rate:11m:omega} \rho_{k+1}-\lambda_1&\le \omega(\sigma,\wtd C)(\rho_k-\lambda_1) \quad\text{if the equality in \cref{eq: thm:convergence-rate:11m} holds},
	\end{align}
	\end{subequations}
	where $\chi(\cdot,\cdot),\omega(\cdot,\cdot)$ are as in \cref{eq:chi-psi} and $ \wtd C\le \Delta^2+O(\rho_k-\lambda_1), \sigma=(1+\abs{\gamma_2}+\dots+\abs{\gamma_{m_h}})^2 .$
 Here $\gamma_j=\frac{(r_{{k-j}}-r_{{k-j+1}})^{\HH}r_k }{r_k^{\HH}r_k }$ for $j=1,\dots,m_h$.
\end{theorem}

	The parameter $\sigma$ relies on the historical terms that is highly related to the choice of the initial $x_0$ and vary significantly case by case.
	Hence it is very difficult, if not impossible, to derive any bound for $\sigma$ except the trivial one $\sigma>1$.
	We will show several $\sigma$'s in some typical examples in \cref{sec:numerical-experiments}.
More can be said on $\sigma$.
Note 
\begin{align*}
	\chi(\sigma,C)&=\frac{2\sigma+C-2\sqrt{\sigma(\sigma+C)}}{2\sigma+1-2\sqrt{\sigma(\sigma+C)}}
= \frac{C^2}{1+{(1-C)}\left(2\sigma-1+2\sqrt{\sigma(\sigma+C)}\right)},
\\
\omega(\sigma,C)&=\frac{\sigma+C-(1-C)\sqrt{\sigma(\sigma+C)}}{1+\sigma(2-C)}
=\frac{C}{1+(1-C)\sqrt{\frac{\sigma}{\sigma+C}}},
\end{align*}
are both monotonically decreasing with respect to $\sigma$ and monotonically increasing with respect to $C$.
The least possible decrease happens in the case that $r_k\perp \range([r_{k-1},\dots,r_{k-m_h}])$, which is usually expected, as in the practice in CG performing on the linear system $Ax=b$, or that in Lanczos method for Hermitian eigenvalue problem,
because this condition guarantees that the process finishes in finite iterations.
It means that at both cases the convergence behavior is so good that it is useless to adopt in more historical terms.
This implies that \emph{the worse the historical convergence is, the more effective adding more historical terms is}.

	It is expected that $\chi(\sigma,\wtd C)\le \wtd C^{1+m_h}$; that is,  the $(1+m_h)$-step asymptotic convergence rate of $\LOCG(1,1,m_h)$ is smaller than that of SD. 
	However, this is not true for any $\sigma$ except $m_h=1$, which is shown in \cref{sec:a-simple-case}.
	For a small $\kappa$ resulting in $\wtd C\approx 0$, $\chi(\sigma,\wtd C)= O(\wtd C^2)\gg \wtd C^{1+m_h}$,
	which tells the convergence rate is much overestimated.
	Conversely, for a large $\kappa$ resulting in $\wtd C\approx 1$, $\chi(\sigma,\wtd C)\approx 1-\frac{4(\sqrt{\sigma+\wtd C}+\sqrt{\sigma})^2}{\kappa+1}, \wtd C^{1+m_h}\approx 1-\frac{4(1+m_h)}{\kappa+1}$,
	and the comparison is reduced to the relation between $\sigma$ and $m_h$,
	which highly relies on the detailed spectrum of $A$ and the clear relation seems very difficult to be obtained.
	In particular, for the case $m_h\le 4$, 
	$\sigma\ge \frac{m_h^2}{4(1+m_h)}$, and thus $\chi(\sigma,\wtd C)\le \wtd C^{1+m_h}$. 
	Anyway, it is good news that $\omega(\sigma,\wtd C)<\omega(1,\wtd C)<\wtd C$, which tells that the $1$-step asymptotic convergence rate of $\LOCG(1,1,m_h)$ is better than that of $\LOCG(1,1,1)$, not to mention SD.

	On the other hand, the overestimate can be treated by more careful consideration. 
	Note that \cref{eq:optim:m} guarantees $\eta<\wtd C$ and $\eta+\frac{\eta_{-1}(1-\eta)^2}{1-\eta_{-1}}<\wtd C$, which ultimately leads $\eta\eta_{-1}\le \chi(1,\wtd C)$ as in \cref{sec:a-simple-case}.
	Then besides \cref{eq: thm:convergence-rate:11m}, it also holds that
	\begin{align*}
		\rho_{k+1}-\lambda_1 &\le \chi_1(\wtd C)(\rho_{k-m_h}-\lambda_1), \qquad \text{and together}
		\\
		\rho_{k+1}-\lambda_1&\le \min\set{\chi(\sigma,\wtd C), \chi_1(\wtd C)}(\rho_{k-m_h}-\lambda_1).
	\end{align*}
	where 
	\begin{align*}
		\chi_1( C)&=\chi(1, C)^{\lfloor{(1+m_h)/2}\rfloor} C^{(1+m_h)-2\lfloor{(1+m_h)/2}\rfloor}
		\\&= { C^{1+m_h}}\left/{[1+{(1- C)}(1+2\sqrt{1+ C})]^{\lfloor{(1+m_h)/2}\rfloor}}\right.
		.
	\end{align*}
	This gives the expected relation that the $(1+m_h)$-step asymptotic convergence rate $\le \wtd C^{1+m_h}$.

Using the techniques shown in \cref{sec:larger-krylov-subspace,sec:more-historical-terms} above,
we are able to establish a general convergence bound for $\LOCG(1,m_e,m_h)$.
\begin{theorem}\label{thm:convergence-rate:1mm}
	Suppose $\lambda_1\le\rho_0<\lambda_2$.
	For the sequences $\set{\rho_{k}},\set{x_{k}}$ produced by $\LOCG(1,m_e,m_h)$,
	\begin{subequations}
	\begin{align}\label{eq: thm:convergence-rate:1mm}
			\rho_{k+1}-\lambda_1&\le \chi(\sigma,\wtd C)
			(\rho_{k-m_h}-\lambda_1), \quad \text{and}
			\\\label{eq: thm:convergence-rate:1mm:omega} \rho_{k+1}-\lambda_1&\le\omega(\sigma,\wtd C)(\rho_k-\lambda_1) \quad\text{if the equality in \cref{eq: thm:convergence-rate:1mm} holds},
\end{align}
	\end{subequations}
	where $\chi(\cdot,\cdot),\omega(\cdot,\cdot)$ are as in \cref{eq:chi-psi} and $
	\wtd C\le \cheb_{m_e}(\Delta^{-1})^{-2}+O(\rho_k-\lambda_1),
	\sigma=(1+\abs{\wtd \gamma_2}+\dots+\abs{\wtd \gamma_{m_h}})^2,
	\wtd\gamma_j=\frac{(r_{k-j}-r_{k-j+1})^{\HH}q_k}{\N{r_k}^2}
	$, in which $q_k$ is as in \cref{eq:qk}.
\end{theorem}

\subsection{More eigenvalues:  $\LOCG(n_b,1,1)$}\label{sec:more-eigenvalues}
For this case, the search subspace is $\range([X_k,r(X_k),X_{k-1}])$ where and the goal is $n_b$ smallest eigenvalues $\lambda_1,\dots,\lambda_{n_b}$.

In the $k$th iteration,
write $R_k=r(X_k)
,P_k=P(X_k)$, and then the orthogonal basis of the search subspace is $X_k, R_k, (I-P_k)X_{k-1}$ for $R_k^{\HH}X_k=R_k^{\HH}X_{k-1}=0$.
Since only the column spaces of $X_k,X_{k-1},X_{k+1}$ are effective in the Rayleigh quotient, we may rewrite the iteration as 
\[
	\what X_{k+1}=X_k+R_ka_k+(I-P_k)X_{k-1}b_k,
	\quad
S_{k}=(\what X_{k+1}^{\HH}\what X_{k+1})^{1/2},
\quad
	X_{k+1} = \what X_{k+1}S_{k}^{-1}Q_{k}^{\HH},
\]
where $(a_k, b_k)$ minimizes the function $\varphi( a, b)=\rho(X_k+ R_ka + (I-P_k)X_{k-1}b)$,
and $Q_{k}$ is a unitary matrix such that $X_{k+1}^{\HH}X_{k+1}=I_{n_b},X_{k+1}^{\HH}AX_{k+1}=\rho_{k+1}$ is a diagonal matrix whose diagonal entries are nondecreasing.
Let the targets $\Lambda_1=\diag(\lambda_1,\dots,\lambda_{n_b}),U=[u_1,\dots,u_{n_b}]$ where $Au_j=u_j\lambda_j$.
Write $\varepsilon_k=X_k^{\HH}\opL_U(X_k)=\rho_k-\Lambda_1$ and then $\delta_k=\varepsilon_k-\varepsilon_{k+1}$.

In order to make things easy to follow, we will consider the new iterative vectors sequentially:
first consider the first vector $x_{k+1,(1)}$ by 
\begin{equation}\label{eq:sequential:1}
	x_{k+1,(1)}=x_{k,(1)}+\alpha_kr_{k,(1)}+(I-P(x_{k,(1)}))W_{k,(1)}b_{k,(1)},
\end{equation}
where $
	W_{k,(1)}=\left[
		 X_{k-1}, X_{k,(2:n_b)}, (I-P(r_{k,(1)}))R_{k,(2:n_b)}
		 \right]$;
	then consider the second vector $x_{k+1,(2)}$ by 
\begin{equation}\label{eq:sequential:2}
	\what x_{k+2,(2)}=\what x_{k,(2)}+\alpha_k\what r_{k,(2)}+(I-P(\what x_{k,(2)}))\what W_{k,(2)}b_{k,(2)}, 
\end{equation}
where $\what v=(I-P(x_{k+1,(1)}))v$, namely this search is done in the orthogonal complement $\range(x_{k+1,(1)})^{\perp}$; 
then search $x_{k+1,(3)}$ in $\range([x_{k+1,(1)}, \what x_{k+1,(2)}])^{\perp}$ and so forth.
Under the help of minimax principles, this process is essentially equivalent to the original one.

Now focus on \cref{eq:sequential:1}.
By following the same procedure used to derive the asymptotic convergence rate of $\LOCG(1,1,1)$, we obtain that of $\LOCG(n_b,1,1)$.
\begin{theorem}\label{thm:convergence-rate:m11}
	Suppose $\lambda_1\le\rho_{0,(1)}<\lambda_2$.
	For the sequences $\set{\rho_{k,(1)}}$, $\set{x_{k,(1)}}$  produced by $\LOCG(n_b,1,1)$,
	\begin{subequations}
	\begin{align}\label{eq: thm:convergence-rate:m11}
		\rho_{k+1,(1)}-\lambda_1&\le \chi(\sigma_{(1)},\wtd C_{(1)})
			(\rho_{k-1,(1)}-\lambda_1), \quad \text{and}
			\\\label{eq: thm:convergence-rate:m11:omega} \rho_{k+1,(1)}-\lambda_1&\le \omega(\sigma_{(1)},\wtd C_{(1)})(\rho_{k,(1)}-\lambda_1) \quad\text{if the equality in \cref{eq: thm:convergence-rate:m11} holds},
	\end{align}
	\end{subequations}
	where $\chi(\cdot,\cdot),\omega(\cdot,\cdot)$ are as in \cref{eq:chi-psi} and $
	\wtd C_{(1)}\le \Delta_{(1)}^2+O(\rho_{k,(1)}-\lambda_1),
	\sigma_{(1)}=1+\gamma_{k,(1)}^2
	,$ in which $\gamma_{k,(1)}$ is defined as in \cref{eq:gamma:m11} or \cref{eq:gamma:m11:shrink}.
\end{theorem}
	The parameter $\gamma_{k,(1)}$ depends on many terms.
Let 
\begin{align*}
	E_k&=\rho_k-\rho_{k+1,(1)}I=\begin{bmatrix}
	\delta_{k,(1)}&\\&E_{k,(2:n_b)}
\end{bmatrix},
\\
\Gamma_k&=
X_k^{\HH}X_{k-1}=\begin{bmatrix}
	x_{k,(1)}^{\HH}X_{k-1} \\
	X_{k,(2:n_b)}^{\HH}X_{k-1} 
\end{bmatrix}=:\begin{bmatrix}
g_k^{\HH}\\
\Gamma_{k,(2:n_b)}
\end{bmatrix}=:\begin{bmatrix}
 \ol g_{k,(1)} & g_{k,(2:n_b)}^{\HH}\\
 h_k  &  H_{k,(2:n_b)}
\end{bmatrix},\qquad g_{k,(1)}\in \mathbb{C}.
\end{align*}
Then
	\begin{align}
	\gamma_{k,(1)}^2
	&=
	\delta_{k-1,(1)}\left(\left[\frac{h_k}{\delta_{k-1,(1)}}+\what\Phi_1\right]^{\HH}\left[M_k-\frac{h_kh_k^{\HH}}{\delta_{k-1,(1)}}-\what\Psi_{21}\what\Psi_{22}\what\Psi_{21}^{\HH}\right]^{-1}\left[\frac{h_k}{\delta_{k-1,(1)}}+\what\Phi_1\right]+\what\Phi_2^{\HH}\what\Psi_{22}\what\Phi_2\right)
	\nonumber\\&=O(1)
	,\label{eq:gamma:m11}
	\end{align}
where
\begin{subequations}\label{eq:gamma:terms}
\begin{align}
	M_k&=E_{k,(2:n_b)}^{-1}-H_{k,(2:n_b)}E_{k,(2:n_b)}^{-1}H_{k,(2:n_b)}^{\HH}+E_{k,(2:n_b)}^{-1}\delta_{k-1,(2:n_b)}[E_{k,(2:n_b)}+\delta_{k-1,(2:n_b)}]^{-1},\\
	\what\Psi_{22}&=R_{\perp}(R_{\perp}^{\HH}[A-\rho_{k+1,(1)}I]R_{\perp})^{-1}R_{\perp}^{\HH}=O(1), \qquad R_{\perp}=(I-P(r_{k,(1)}))R_{k,(2:n_b)},\\
	\what\Psi_{21}&= E_{k,(2:n_b)}^{-1}R_{\perp}^{\HH}=O(\N{\varepsilon_k^{1/2}}),\\
	\what\Phi_1&=
		\frac{1}{r_{k,(1)}^{\HH}r_{k,(1)}}E_{k,(2:n_b)}^{-1}\left[R_{k,(2:n_b)}^{\HH}-R_{\perp}^{\HH}\what\Psi_{22}(A-\rho_{k,(1)}I)\right]r_{k,(1)}
	=O(1)
	,
	\\
	\what\Phi_2
	&=\frac{1}{r_{k,(1)}^{\HH}r_{k,(1)}}{(A-\rho_{k,(1)}I)r_{k,(1)}}
	=\frac{O(1)}{\delta_{k,(1)}^{1/2}}
	.
\end{align}
\end{subequations}

\begin{remark}\label{rk:h:new-order}
	As a by-product, a new sharp estimate of the difference between the approximations of two successive steps is also established.
For $i=1,\dots,n_b$, write $x_{k,(i)}=u_i\xi_{k,(i)}+v_i\varpi_{k,(i)}+w_i\zeta_{k,(i)}$, where $v_i\in\range(U),v_i\perp u_i, w_i\perp \range(U)$.
Previously,
\begin{enumerate}
	\item In the general case, 
		$1-\abs{\xi_{k,(i)}}^2=\abs{\varpi_{k,(i)}}^2+\abs{\zeta_{k,(i)}}^2=O(\varepsilon_{k,(i)})$. 

	\item In the ``quite near'' case that $\N{\varepsilon_k}$ is sufficiently small, according to \cite[Lemma~5]{longsineM1980simultaneous}, 
		$\abs{\varpi_{k,(i)}}=O(\abs{\zeta_{k,(i)}}\abs{\zeta_{k,(t)}})=O(\varepsilon_{k,(i)}^{1/2}\N{\varepsilon_{k-1}^{1/2}})$ for some $t\ne i$.
\end{enumerate}
Thus
$x_{k,(i)}^{\HH}x_{k-1,(i)}=1-O(\varepsilon_{k-1,(i)}^{1/2}\varepsilon_{k-1,(i)}^{1/2})$;
for $i\ne j$, $x_{k,(i)}^{\HH}x_{k-1,(j)}=O(\varepsilon_{k,(i)}^{1/2}+\varepsilon_{k-1,(j)}^{1/2})\N{\varepsilon_{k-1}^{1/2}}^{\theta}$ where $\theta=0$ for the general case or $\theta=1$ for the quite near case.
Similarly, $r_{k,(i)}^{\HH}r_{k,(j)}=O(\varepsilon_{k,(i)}^{1/2}\varepsilon_{k,(j)}^{1/2})$;
$r_{k,(i)}^{\HH}(A-\rho_{k,(1)}I)r_{(j)}=O(\varepsilon_{k,(i)}^{1/2}\varepsilon_{k,(j)}^{1/2})$.
With $\delta_{k,(1)}/\delta_{k-1,(1)},\varepsilon_{k,(1)}/\delta_{k,(1)}$ bounded under mild condition (as in \cref{sec:a-simple-case}),
we have 
\begin{gather}
	g_{k,(1)}=x_-^{\HH}x=1-O(\varepsilon_{k,(1)}),\qquad
	g_{k,(2:n_b)}=X_{k-1,(2:n_b)}^{\HH}x=O(\N{\varepsilon_{k-1}^{(1+\theta)/2}}),
	\nonumber
	\\
	\label{eq:order-h:original}
	h_k =X_{k,(2:n_b)}^{\HH}x_{k-1}=O(\N{\varepsilon_{k-1}^{(1+\theta)/2}}),\qquad
	H_{(2:n_b)} =X_{k,(2:n_b)}^{\HH}X_{k-1,(2:n_b)}=I-O(\N{\varepsilon_{k-1}^{(1+\theta)/2}}).
\end{gather}
Let $\tau_k^2=\frac{h_k^{\HH}M_k^{-1}h_k}{\delta_{k-1,(1)}}$.
In the proof of \cref{thm:convergence-rate:m11}, we will see $0<\tau_k^2<1,M_k\succ0$, which guarantees
\begin{equation}\label{eq:h:new-order}
	h_k=O(\delta_{k-1,(1)}^{1/2}\N{\varepsilon_{k-1}^{(1+\theta)/4}}).
\end{equation}
This is stronger than \cref{eq:order-h:original}! 

On the other hand,
in particular, for $\theta=0$, 
\begin{equation}\label{eq:gamma:m11:shrink}
	\gamma_{k,(1)}^2= \frac{\tau_k^2}{1-\tau_k^2}+\delta_{k-1,(1)}\what\Phi_2^{\HH}\what\Psi_{22}\what\Phi_2
=O(1),
\end{equation}
which is much simpler than \cref{eq:gamma:m11}.
\end{remark}

Return to \cref{thm:convergence-rate:m11}.
Similarly we discuss a little more on $\sigma$.
	Again the parameter $\sigma$ relies on the historical terms highly related to the initial $x_0$ and vary significantly case by case.
	We will show several $\sigma$'s in some typical examples in \cref{sec:numerical-experiments}.
Noticing the monotonicity of the rate as a function of $\sigma$,
the larger $\gamma_{k,(1)}$ is, 
the larger $\sigma$ is, and the smaller the rate is.
This tells that \emph{the inaccurate approximation of the other eigenvalues helps improve the accuracy of the smallest eigenvalue},
which is reasonable, since the inaccurate approximation has large component in the direction of the eigenvector associated with the smallest eigenvalue, and hence helps for its appearance in the search subspace.
This also implies that \emph{usually the smallest eigenvalue is first well-approximated and locked before any other eigenvalues}, as in the practical computation experience.

	Moreover, since $\chi(\sigma, C)\le \chi(1,C)< \left(\frac{C}{2-C}\right)^2 < C^2$,
	and $\omega(\sigma,C)\le \omega(1,C)<C$,
	the asymptotic convergence rate of the smallest eigenvalue by $\LOCG(n_b,1,1)$ is smaller than that of $\LOCG(1,1,1)$ as expected.

For \cref{eq:sequential:2} and more terms, 
using \cref{thm:convergence-rate:m11} produces,
for instance,
\[
	\rho_{k+1,(2)}-\mu_2\le \chi(\sigma_{(2);\mu},\wtd C_{(2);\mu})
	(\rho_{k-1,(2)}-\mu_2),
\]
where $\mu_2$ is the smallest eigenvalue of $(I-x_{k,(1)}x_{k,(1)}^{\HH})A(I-x_{k,(1)}x_{k,(1)}^{\HH})$,
and $\wtd C_{(2);\mu}\le \Delta_{(2)}^2+O(\rho_{k,(2)}-\lambda_1), \sigma_{(2);\mu}=1+\gamma_{k,(2)}^2$, $\gamma_{k,(2)}$ relies on $\delta_{k-1,(3)},\dots,\delta_{k-1,(n_b)}$.
Since $0\le \lambda_2-\mu_2\le \varepsilon_{k,(1)}$,
the convergence for the second smallest eigenvalue relies on that of the first one, which is usually detected in the practical computation experience.
Moreover, since $\frac{\rho_{k+1,(2)}-\mu}{\rho_{k-1,(2)}-\mu}$ is monotonically decreasing with respect to $\mu$,
we have
\[
	\rho_{k+1,(2)}-\lambda_2\le \chi(\sigma_{(2)},\wtd C_{(2)})
	(\rho_{k-1,(2)}-\lambda_2),
\]
which suggests a total estimate below.
\begin{theorem}\label{thm:convergence-rate:m11:total}
	Suppose $\lambda_1\le\rho_{0,(1)}<\lambda_2\le\rho_{0,(2)}<\dots<\lambda_{n_b}\le\rho_{0,(n_b)}<\lambda_{n_b+1}$.
	For the sequences $\set{\rho_{k}}$, $\set{X_{k}}$  produced by $\LOCG(n_b,1,1)$,
	\begin{align*}
		\tr(\rho_{k+1}-\Lambda_1)&\le \chi(\sigma,\wtd C)\tr(\rho_{k-1}-\Lambda_1).
	\end{align*}
	where $\chi(\cdot,\cdot),\omega(\cdot,\cdot)$ are as in \cref{eq:chi-psi} and $
	\wtd C\le \max_i\Delta_{(i)}^2+O(\N{\varepsilon_{k}}),
	\sigma=1+\max_i\gamma_{k,(i)}^2
	.$
\end{theorem}

\section{Numerical Experiments}\label{sec:numerical-experiments}
We present numerical results from three representative examples to illustrate the behavior of LOBPeCG and verify the obtained asymptotic convergence rates together with their tightness.
All experiments were conducted in MATLAB 2023b under the Windows 11 Professional 64-bit operating system on a PC with an Intel Core i7-11370H processor at 3.30GHz and 32GB RAM. 

The three test problems are denoted as Laplacian 50, Cluster-Outlier, and Outlier-Cluster.

For each example, we test $13$ different triples $(n_b,m_e,m_h)$:
\begin{enumerate}
	\item $\LOCG(1,m_e,m_h)$ for $m_e=1,2,3,m_h=0,1,2$: to recognize the effect of increasing $m_e$ and/or $m_h$;
	\item $\LOCG(n_b,1,m_h)$ for $n_b=1,2,3,m_h=0,1$: to recognize the effect of increasing $n_b$ for BSD and LOBCG.
\end{enumerate}
Note that the two classes share two common triples.

In each example, $20$ independent trials were run with different initial $X_0$'s and present only $3$ typical ones out of them to save space. In detail, in each trial $X_0$ is constructed by taking the first $n_b$ columns of an orthonormal matrix $\wtd X_0$ of size $n\times 3$, uniformly sampled from the Stiefel manifold.

The iteration stops when either of the following criteria is met: the decrease of Rayleigh quotient $\rho_{k-1}-\rho_k < 10^{-15}\abs{\rho_k}$, or nonreal Ritz values are encountered during computation.

\begin{figure}[hp]
	\centering
	\includegraphics[height=0.325\textheight]{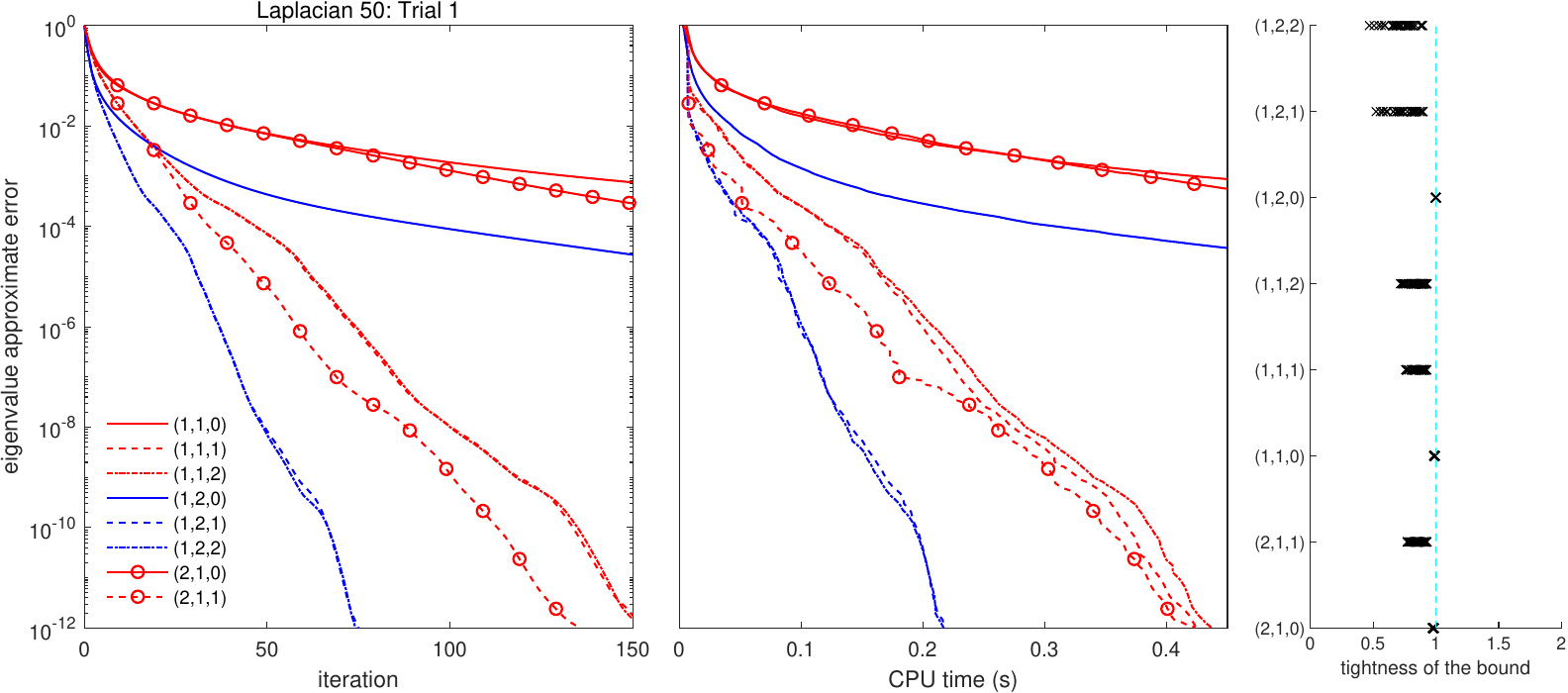}

\smallskip

	\includegraphics[height=0.325\textheight]{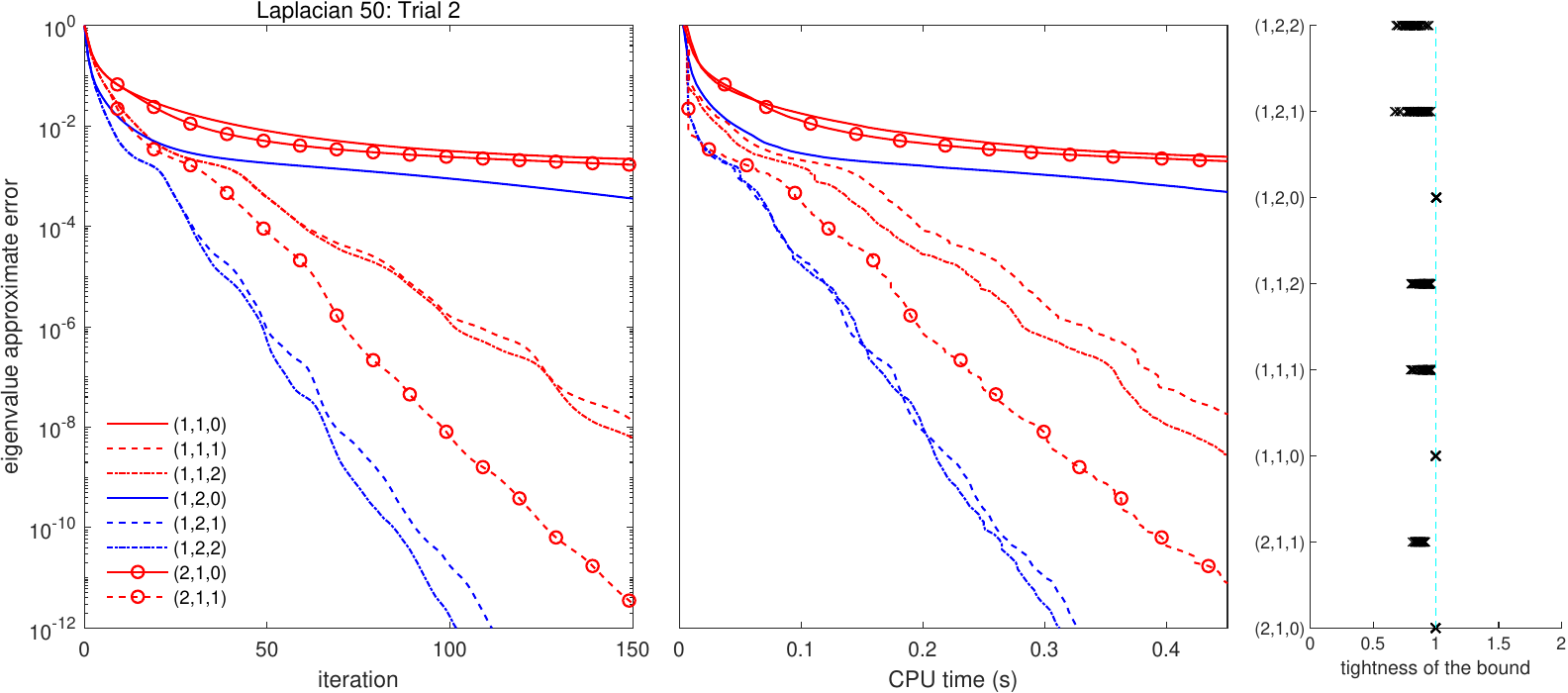}

\smallskip

	\includegraphics[height=0.325\textheight]{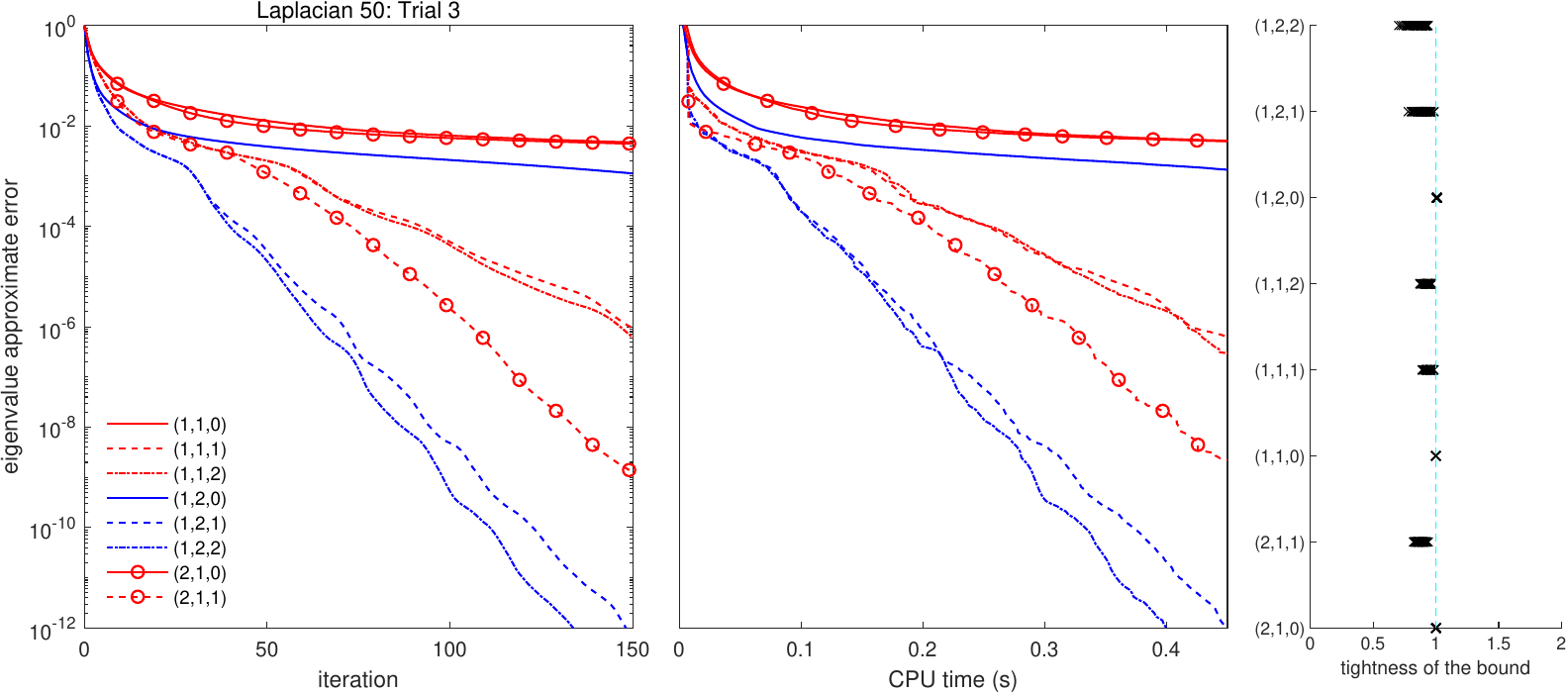}
	\caption{Laplacian 50}
	\label{fig::laplacian50-1}
\end{figure}

\begin{figure}[hp]
	\centering
	\includegraphics[height=0.325\textheight]{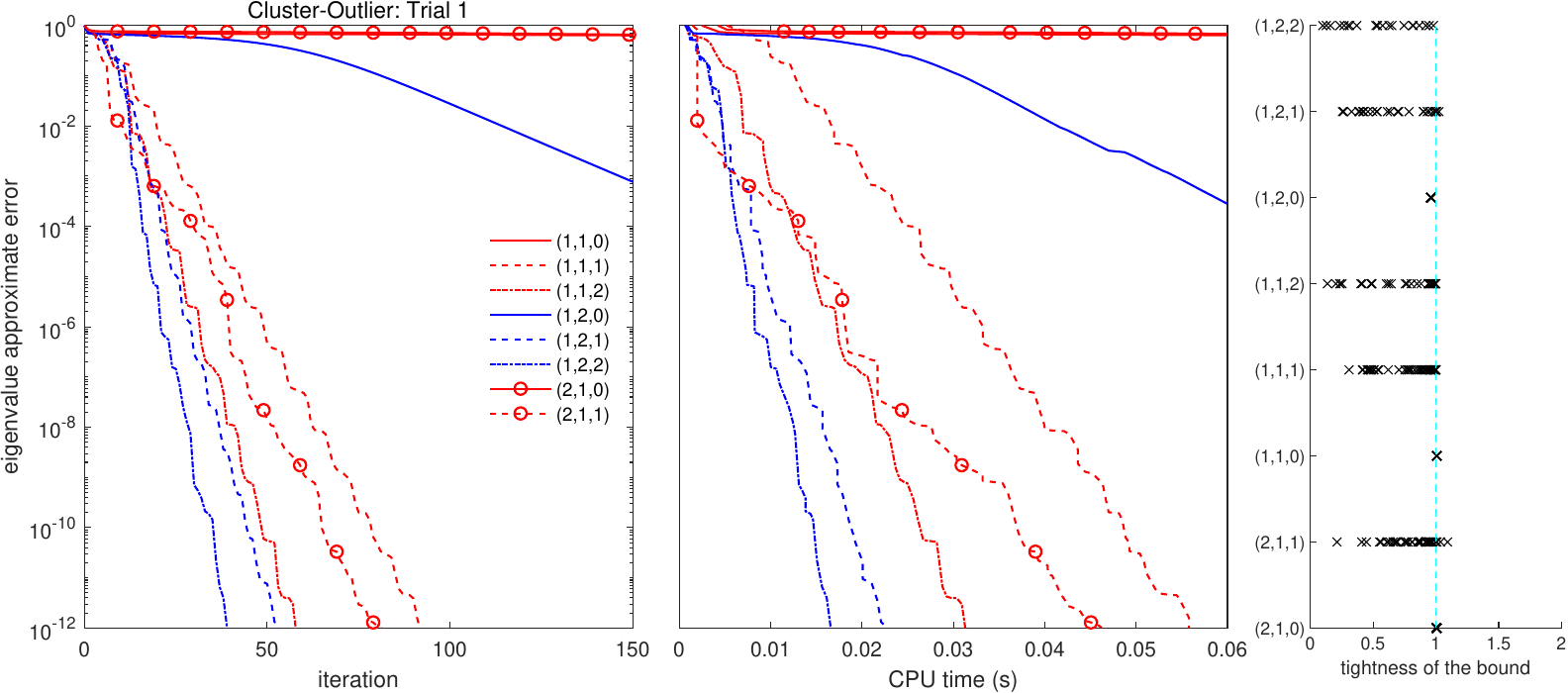}

\smallskip

	\includegraphics[height=0.325\textheight]{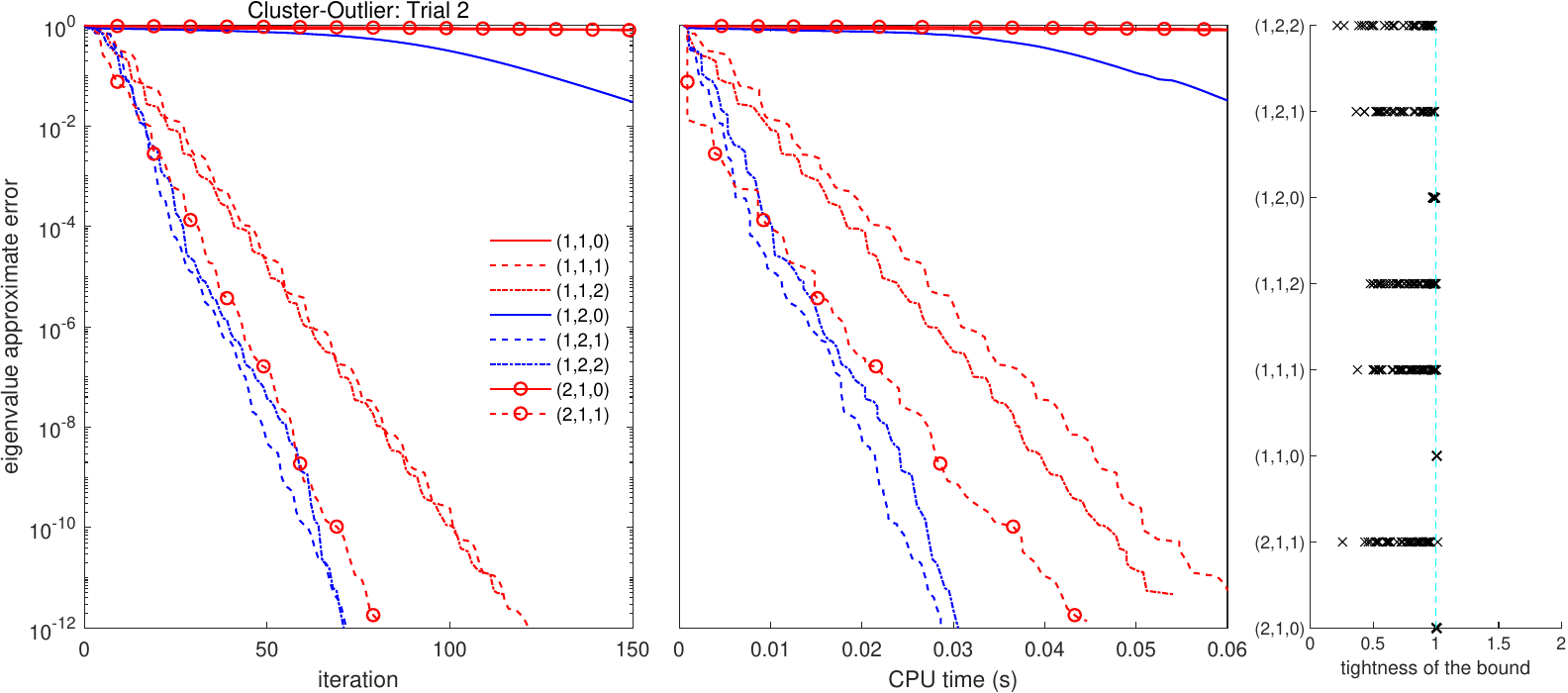}

\smallskip

	\includegraphics[height=0.325\textheight]{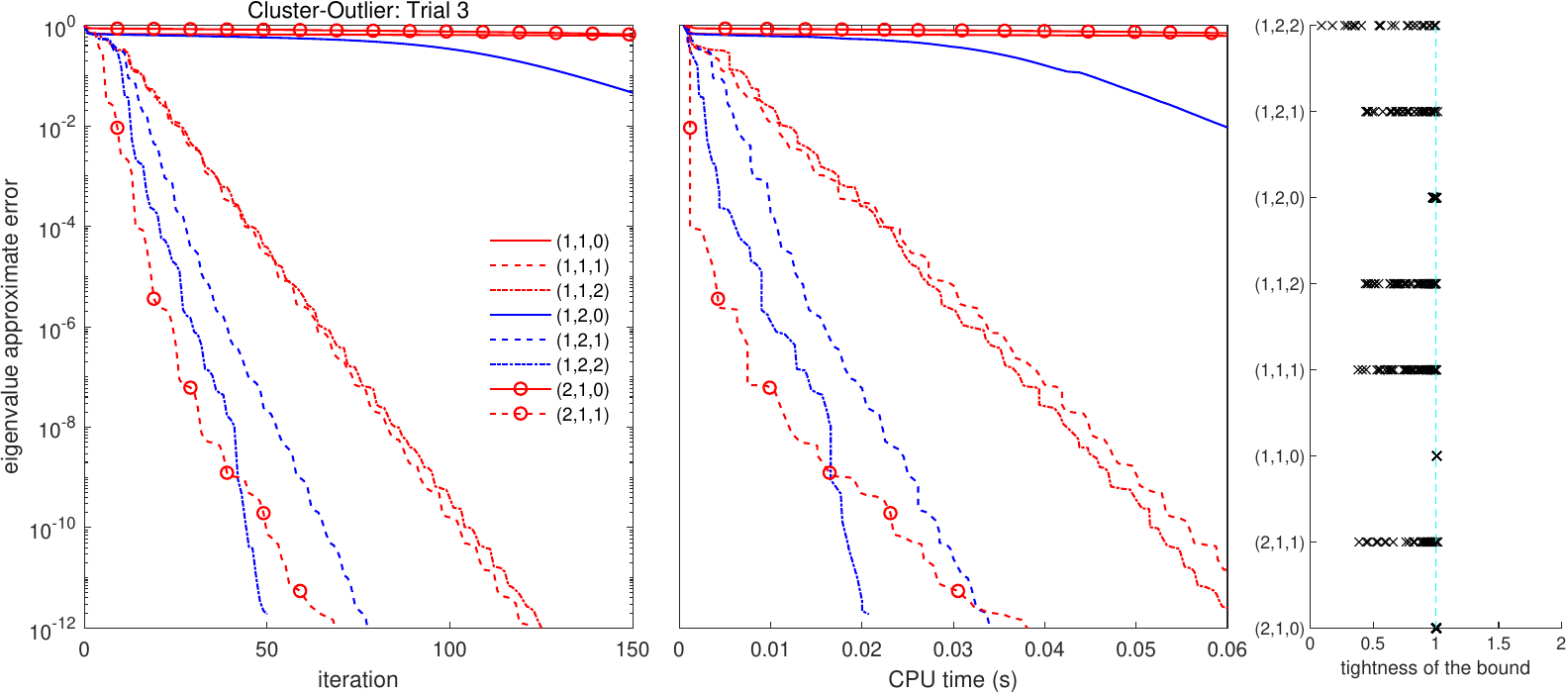}
	\caption{Cluster-Outlier}
	\label{fig::worst-1}
\end{figure}

\begin{figure}[hp]
	\centering
	\includegraphics[height=0.325\textheight]{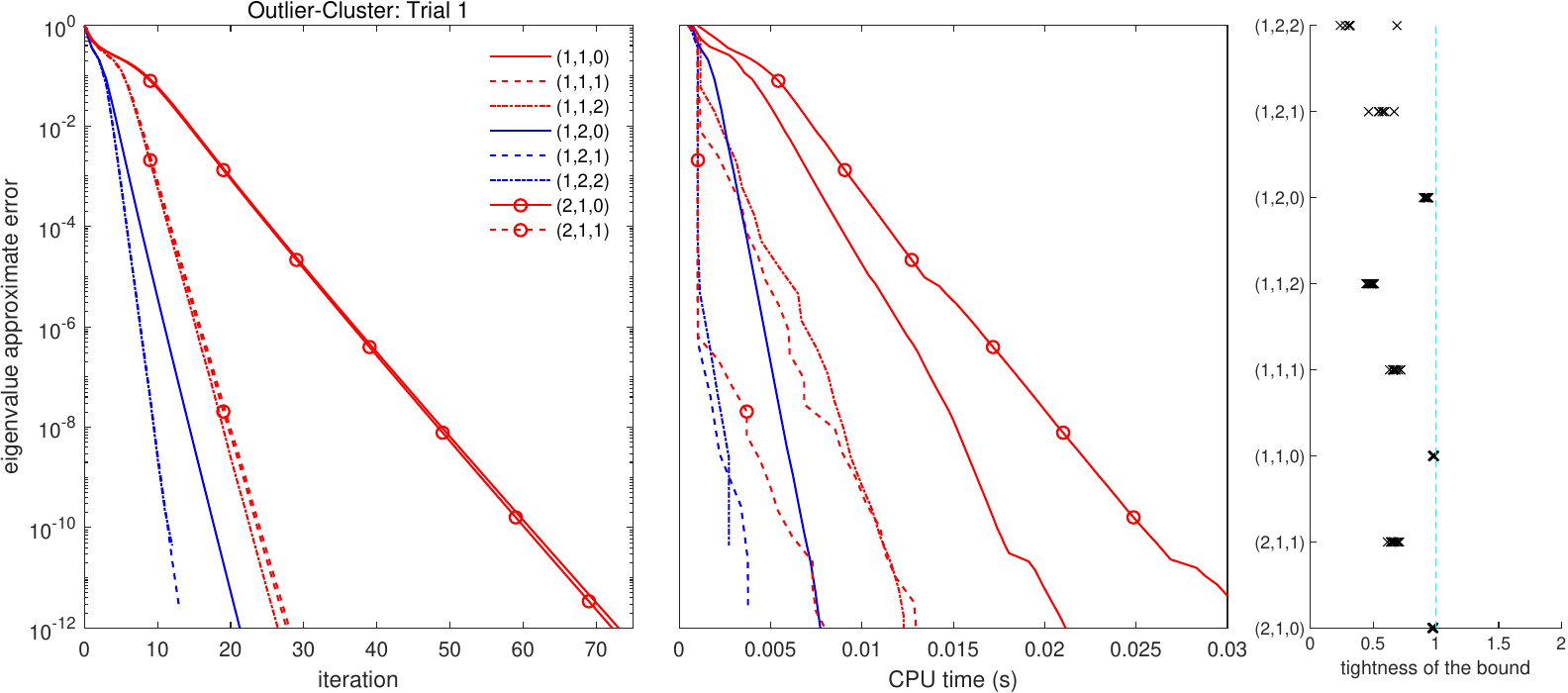}

\smallskip

	\includegraphics[height=0.325\textheight]{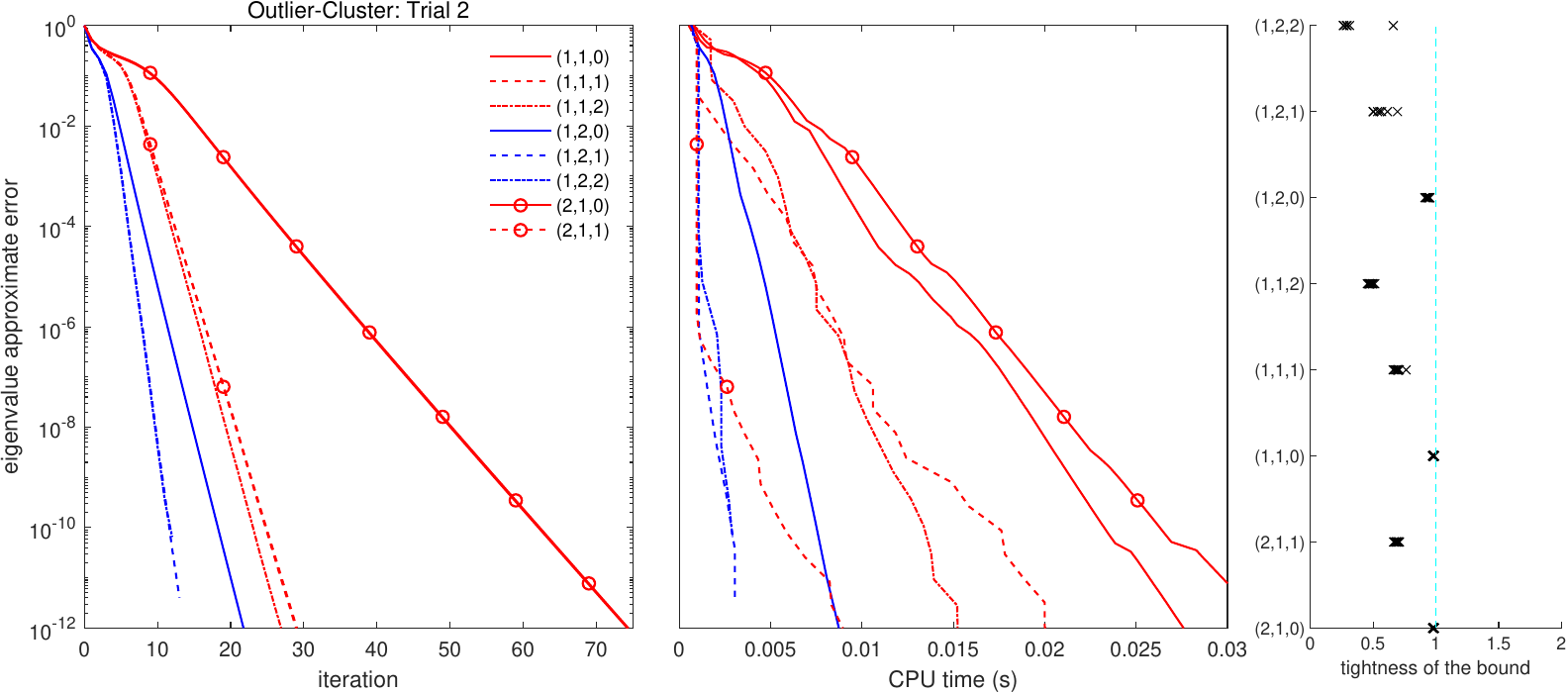}

\smallskip

	\includegraphics[height=0.325\textheight]{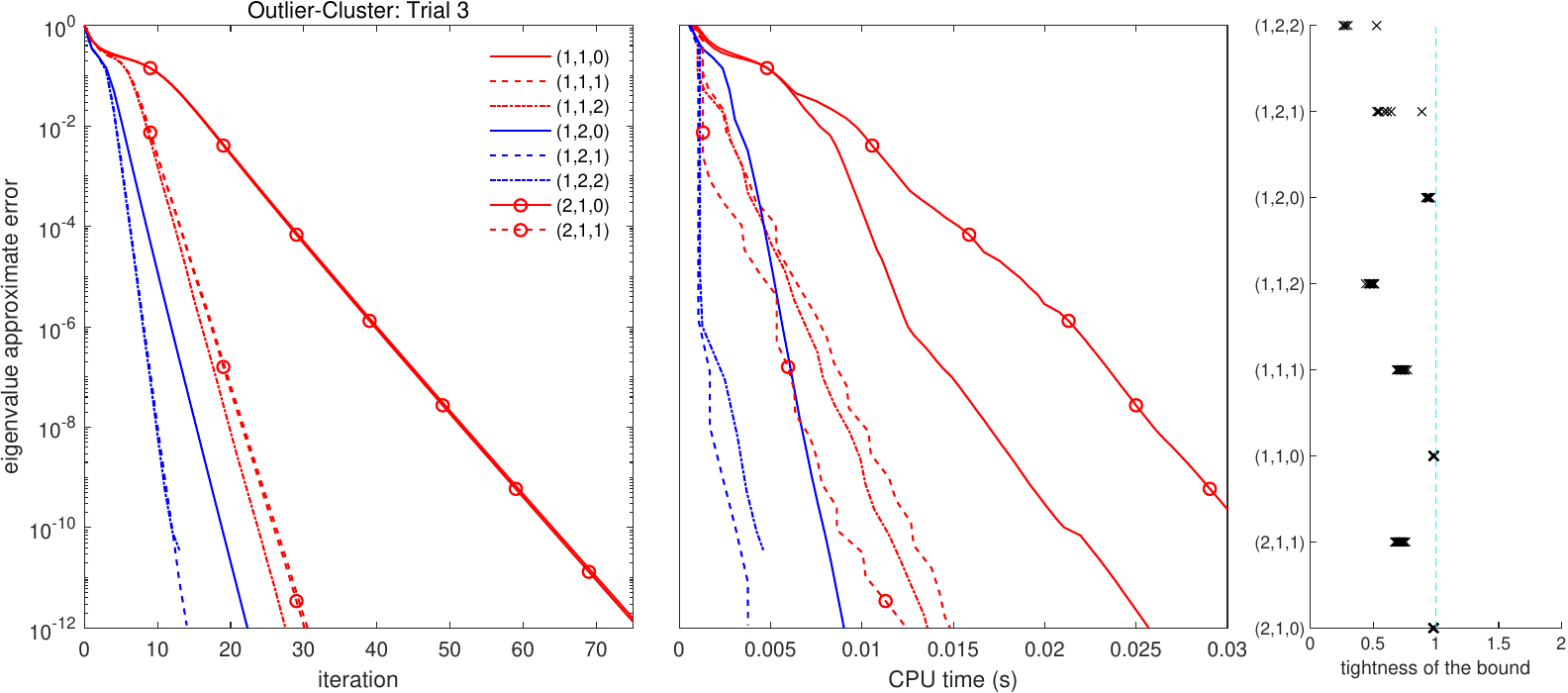}
	\caption{Outlier-Cluster}
	\label{fig::best-1}
\end{figure}

The convergence behavior of each example are given in \cref{fig::laplacian50-1,fig::worst-1,fig::best-1} correspondingly.
Each figure consists of $3$ subfigures for the $3$ trials, of which each has three parts:
the left part illustrates the performance for iteration vs.\ the relative eigenvalue approximate error $\frac{\rho_k-\lambda_1}{\rho_0-\lambda_1}$;
the middle part illustrates the performance for CPU time (in seconds) vs.\ the relative eigenvalue approximate error;
the right part demonstrates the validity and tightness of the obtained asymptotic convergence rate by the quantities $\frac{\rho_k-\lambda_1}{\sqrt{\chi(\sigma,C)}(\rho_{k-1}-\lambda_1)}$ in the last 50\% iterations.
Since the actual $C$ or $\wtd C$ is not in hand, we have to use the $C$ generated by Chebyshev uniform bound instead.
Note that for the triple $(n_b,1,m_h)$, only the error of the smallest eigenvalue $\lambda_1$ is reported.
To make figures readable, only $8$ out of $13$ triples are shown in the figures.

\begin{example}[Laplacian 50]\label{eg:laplacian-50}
	The example is a classic test problem originated from a finite-difference discretization of the Laplacian operator on a unit square with an $N\times N$ uniform grid. 
	In this example, $A$ is a tridiagonal matrix with diagonal entry $-4$ and sub/super-diagonal entry $1$.
	Its eigenvalues are expressed by
	$-4+2\cos\frac{i\pi}{N+1}+2\cos\frac{j\pi}{N+1}$ for $i,j=1,\dots,N$.
	Here $N=50, n=N^2=2500$ and LOBPeCG is tested to compute the smallest eigenvalue.
	$\kappa\doteq702.5410,
	\Delta \doteq.9972,
	C_{\cheb}=\cheb_{m_e}(\Delta^{-1})^{-2}\doteq.9943,.9775,.9504$ for $m_e=1,2,3$.
	The results are shown in \cref{fig::laplacian50-1}.
\end{example}
\begin{example}[Cluster-Outlier]\label{eg:cluster-outlier}
	The example is artificially designed for a complicated behavior.
	In this example, $A=Q\Lambda Q^H$, where $\Lambda=\diag(-1,0,0.001,\dots,0.994, 2^6,2^7,\dots,2^9)$ and $Q$ is an orthogonal matrix uniformly sampled in the special orthogonal group. 
	Here $n=1000$ and LOBPeCG is tested to compute the smallest eigenvalue $-1$.
	$\kappa=513,
	\Delta \doteq.9961,
	C_{\cheb}=\cheb_{m_e}(\Delta^{-1})^{-2}\doteq.9922,.9694,.9329$ for $m_e=1,2,3$.
	The results are shown in \cref{fig::worst-1}.
\end{example}
\begin{example}[Outlier-Cluster]\label{eg:outlier-cluster}
	The example is artificially constructed for small $\kappa$ and $C$.
	In this example, $A$ is constructed as that in \cref{eg:outlier-cluster} except $\Lambda=\diag(1,1.1,1.101,\dots,2.098)$.
	Here $n=1000$ and LOBPeCG is tested to compute the smallest eigenvalue $1$.
	$\kappa=10.98,
	\Delta \doteq.8331,
	C_{\cheb}=\cheb_{m_e}(\Delta^{-1})^{-2}\doteq.6940,.2824,.0908$ for $m_e=1,2,3$.
	The results are shown in \cref{fig::best-1}.
\end{example}

For the performance of all the variants of LOBPeCG,
we can observe that
\begin{enumerate}
\item Increasing $m_h$ consistently improves convergence, as it introduces little additional computational cost while potentially offering significant convergence acceleration.
\item Increasing $m_e$ and/or $n_b$ generally improves convergence, but the associated increase in computational cost must be weighed against the benefit.
\item The relative efficiency of increasing $m_e$ versus $n_b$ is problem-dependent; neither strategy is uniformly superior.
\item Increasing $n_b$ allows the algorithm to compute additional eigenvalues, which is an inherent advantage of the block approach.
\end{enumerate}

For the validity and tightness of the obtained asymptotic convergence rate,
we can see that
our estimate for the upper bound of the rate in $\LOCG(n_b,m_e,m_h)$ is usually good, although sometimes underestimated. 
Despite the computational error, this occasional underestimation can be attributed by the fact that
the comparison is between the actual rate and the averaged rate in the two-step rate $\sqrt{\chi(\sigma,C)}$,
and the underestimation would be averaged by the overestimation for the next/last step.
On the other hand, the actual rate can be much smaller than the estimate, because we have to use $C$ generated from the Chebyshev uniform bound other than the actual $C$ or $\wtd C$ that would be much smaller.

To examine the effect of different $\sigma$'s, \cref{tab:different-sigma-s-in-the-examples} summarizes the min/mean/max values of the quantities $\sigma-1$ in the last 50\% iterations of the $3\times 3$ trials for the following five parameter triples $(1,m_e,2),(n_b,1,1)$ for $m_e=1,2,3,n_b=2,3$.
Any value less than $10^{-8}$ is presented as $0$, while any one more than $1$ is bolded.
The $\sigma$'s of the other triples are omitted because theoretically $\sigma=0$ for $(1,m_e,0)$, $\sigma=1$ for $(1,m_e,1)$.

The data indicate that $\sigma$ is often close to $1$, but for certain cases $\sigma$ is much larger than $1$, leading to much faster convergence, especially in the experiment on Cluster-Outlier $\LOCG(n_b,1,1)$ for $n_b=2,3$.

\begin{table}[hb]
	\centering
	\scriptsize
\begin{tabular}{|@{\,}c@{\,}|@{\,}c@{\,}|*{5}{@{\,}c@{\,}|}}
\hline
\multicolumn{2}{|@{}c@{\!\!}|@{\,}}{$\sigma-1$}& {$(1,1,2)$} &  {$(1,2,2)$}&   {$(1,3,2)$}  &  {$(2,1,1)$} & {$(3,1,1)$}  \\
\hline
\multirow{1}{*}{\rotatebox{90}{~L.~50~}}&Trial 1&
\hmz,\hmz,\hmz&
\hmz,6.8e-7,6.8e-6&
\hmz,2.3e-5,2.6e-4&
\hmz,1.6e-4,1.2e-3&
2.3e-7,2.6e-4,2.1e-3\\
&Trial 2&
\hmz,\hmz,\hmz&
\hmz,\hmz,\hmz&
\hmz,9.9e-8,1.5e-6&
1.8e-8,1.6e-4,8.2e-4&
3.8e-6,3.3e-4,3.9e-3\\
&Trial 3&
\hmz,\hmz,\hmz&
\hmz,\hmz,\hmz&
\hmz,\hmz,\hmz&
\hmz,2.5e-4,4.2e-3&
1.9e-6,2.9e-4,1.8e-3\\
\hline
\multirow{1}{*}{\rotatebox{90}{~C.-O.~}}&Trial 1&
\hmz,\hmz,\hmz&
\hmz,\hmz,\hmz&
\hmz,1.6e-5,2.0e-4&
7.8e-6,3.5e-1,{\bf 2.7e0}&
5.7e-5,8.8e-1,{\bf 4.3e0}\\
&Trial 2&
\hmz,\hmz,\hmz&
\hmz,\hmz,\hmz&
\hmz,5.1e-8,5.1e-7&
4.0e-4,5.8e-1,{\bf 1.5e1}&
1.6e-3,{\bf 3.6e0},{\bf 2.7e1}\\
&Trial 3&
\hmz,\hmz,\hmz&
\hmz,\hmz,1.1e-7&
3.7e-6,6.9e-4,4.2e-3&
5.1e-6,4.7e-1,{\bf 4.0e0}&
7.3e-4,{\bf 1.5e0},{\bf 1.1e1}\\
\hline
\multirow{1}{*}{\rotatebox{90}{~O.-C.~}}&Trial 1&
\hmz,\hmz,\hmz&
\hmz,2.2e-7,1.4e-6&
3.5e-8,3.4e-3,1.3e-2&
\hmz,2.3e-4,1.5e-3&
2.7e-5,7.3e-4,5.9e-3\\
&Trial 2&
\hmz,\hmz,\hmz&
\hmz,4.5e-6,2.9e-5&
3.5e-8,4.0e-6,1.6e-5&
1.5e-7,1.6e-4,6.7e-4&
2.4e-5,3.7e-4,1.6e-3\\
&Trial 3&
\hmz,\hmz,\hmz&
\hmz,1.5e-1,{\bf 1.4e0}&
2.2e-8,1.2e-5,3.6e-5&
2.5e-7,1.1e-4,7.9e-4&
4.2e-5,3.1e-4,7.2e-4\\
\hline
\end{tabular}
	\caption{Different $\sigma$'s in the Examples}
	\label{tab:different-sigma-s-in-the-examples}
\end{table}


\section{Concluding remarks}\label{sec:preconditioners}
In this paper we establish a comprehensive and clear analysis of the LOBPCG method for computing the smallest eigenvalues of a Hermitian matrix.
Clearly all the obtained results are valid for the largest eigenvalues when we carefully deal with the signs of differences between eigenvalues and their approximations.
After analyzing the basic case for positive definite Hermitian matrix $F(\lambda)=A-\lambda I$, we briefly explain the performance of LOBPeCG on similar problems, such as positive definite matrix pair $F(\lambda)=A-\lambda B$ and hyperbolic quadratic Hermitian matrix polynomial $F(\lambda)=\lambda^2A+\lambda B+C$.
Since the extreme eigenvalues of these kinds of problems admit the variational characterization, the results still hold for them when proper shifts and metrics (equivalently the inner product) are used.

Finally on the effect of the preconditioners,
for positive definite Hermitian $K_k$, the $i$th iteration is just equivalent to
the $i$th iteration of the algorithm applied to $K_k^{1/2}F(\lambda)K_k^{1/2}$ without a preconditioner.
Thus the obtained convergence rates above all hold as long as $F(\lambda)=A-\lambda I$ is replaced with $K_k^{1/2}F(\lambda)K_k^{1/2}=K_k^{1/2}AK_k^{1/2}-\lambda K_k$, including \cref{thm:convergence-rate:111,thm:convergence-rate:1m1,thm:convergence-rate:11m,thm:convergence-rate:1mm,thm:convergence-rate:m11,thm:convergence-rate:m11:total}, which can be easily examined.

Briefly speaking, we build a bridge that \emph{quantitatively} links the error of iterates in LOBPeCG with those in polynomial iterations (SD, Chebyshev, etc.), and then give an exact estimate based on the convergence rate relying on Chebyshev polynomials.
This techniques can also be used to analyze those methods proposed by adding a search direction for classic methods to achieve a global quasi-optimality, like methods in \cite{wuXS2019trplk,zhangX2024chebyshev}.
Moreover,
it
should be emphasized that this technique can probably be used to analyze other gradient-type optimization methods, including Nesterov accelerated gradient method and stochastic gradient descent method.

\subparagraph{Acknowledgement}
During the preparation of this work the authors used ChatGPT in order to improve readability and language. After using this service, the authors reviewed and edited the content as needed and take full responsibility for the content of the publication.

The authors thank Prof.\ Zhaojun Bai for the helpful discussions during ICBS 2025.



\bibliographystyle{siamplain}
\bibliography{NLAA/cgconv2-NLAA-min}

\begin{thebibliography}{10}

\bibitem{argentatiKNOZ2015convergence}
{\sc M.~E. Argentati, A.~V. Knyazev, K.~Neymeyr, E.~E. Ovtchinnikov, and
  M.~Zhou}, {\em Convergence theory for preconditioned eigenvalue solvers in a
  nutshell}, Foundations of Computational Mathematics, 17 (2017), pp.~713--727.

\bibitem{bennerL2022convergence}
{\sc P.~Benner and X.~Liang}, {\em Convergence analysis of vector extended
  locally optimal block preconditioned extended conjugate gradient method for
  computing extreme eigenvalues}, Numer. Lin. Alg. Appl., 29 (2022), p.~e2445,
  \url{https://doi.org/10.1002/nla.2445}.
\newblock 24 pages.

\bibitem{bergamaschiGP1997asymptotic}
{\sc L.~Bergamaschi, G.~Gambolati, and G.~Pini}, {\em Asymptotic convergence of
  conjugate gradient methods for the partial symmetric eigenproblem}, Numer.
  Lin. Alg. Appl., 4 (1997), pp.~69--126.

\bibitem{bergamaschiP2002numerical}
{\sc L.~Bergamaschi and M.~Putti}, {\em Numerical comparison of iterative
  eigensolvers for large sparse symmetric matrices}, Comput. Methods Appl.
  Mech. Engrg., 191 (2002), pp.~5233--5247.

\bibitem{bradburyF1966new}
{\sc W.~Bradbury and R.~Fletcher}, {\em New iterative methods for solution of
  the eigenproblem}, Numer. Math., 9 (1966), pp.~259--267.

\bibitem{cohen1972rate}
{\sc A.~Cohen}, {\em Rate of convergence of several conjugate gradient
  algorithms}, SIAM J. Numer. Anal., 9 (1972), pp.~248--259.

\bibitem{crowderW1972linear}
{\sc H.~Crowder and P.~Wolfe}, {\em Linear convergence of the conjugate
  gradient method}, IBM Journal of Research and Development, 16 (1972),
  pp.~431--433.

\bibitem{daiY1999nonlinear}
{\sc Y.~Dai and Y.~Yuan}, {\em A nonlinear conjugate gradient method with a
  strong global convergence property}, SIAM J. Optim., 10 (1999), pp.~177--182.

\bibitem{daniel1967conjugate}
{\sc J.~Daniel}, {\em The conjugate gradient method for linear and nonlinear
  operators equations}, SIAM J. Numer. Anal., 4 (1967).

\bibitem{fletcherR1964function}
{\sc R.~Fletcher and C.~Reeves}, {\em Function minimization by conjugate
  gradients}, Comput. J., 7 (1964), pp.~149--154.

\bibitem{gilbertN1992global}
{\sc J.~Gilbert and J.~Nocedal}, {\em Global convergence properties of
  conjugate gradient methods for optimization}, SIAM J. Optim., 2 (1992),
  pp.~21--42.

\bibitem{golubY2002inverse}
{\sc G.~H. Golub and Q.~Ye}, {\em An inverse free preconditioned {K}rylov
  subspace method for symmetric generalized eigenvalue problems}, SIAM J. Sci.
  Comput., 24 (2002), pp.~312--334.

\bibitem{hestenesS1952methods}
{\sc M.~Hestenes and E.~Stiefel}, {\em Methods of conjugate gradients for
  solving linear systems}, J. Res. Nat. Bur. Standards, 49 (1952),
  pp.~409--436.

\bibitem{knyazev2001toward}
{\sc A.~V. Knyazev}, {\em Toward the optimal preconditioned eigensolver:
  Locally optimal block preconditioned conjugate gradient method}, SIAM J. Sci.
  Comput., 23 (2001), pp.~517--541.

\bibitem{knyazevN2003geometricIII}
{\sc A.~V. Knyazev and K.~Neymeyr}, {\em A geometric theory for preconditioned
  inverse iteration {III}: A short and sharp convergence estimate for
  generalized eigenvalue problems}, Linear Algebra Appl., 358 (2003),
  pp.~95--114.

\bibitem{longsineM1980simultaneous}
{\sc D.~E. Longsine and S.~F. McCormick}, {\em Simultaneous {R}ayleigh-quotient
  minimization methods for ${Ax=\lambda Bx}$}, Linear Algebra Appl., 34 (1980),
  pp.~195--234.

\bibitem{neymeyr2001geometricI}
{\sc K.~Neymeyr}, {\em A geometric theory for preconditioned inverse iteration,
  i: Extrema of the {Rayleigh} quotient}, Linear Algebra Appl., 332 (2001),
  pp.~61--85.

\bibitem{neymeyr2012geometric}
{\sc K.~Neymeyr}, {\em A geometric convergence theory for the preconditioned
  steepest descent iteration}, SIAM J. Numer. Anal., 50 (2012), pp.~3188--3207.

\bibitem{neymeyrOZ2011convergence}
{\sc K.~Neymeyr, E.~E. Ovtchinnikov, and M.~Zhou}, {\em Convergence analysis of
  gradient iterations for the symmetric eigenvalue problem}, SIAM J. Matrix
  Anal. Appl., 32 (2011), pp.~443--456.

\bibitem{ovtchinnikov2008jacobiI}
{\sc E.~E. Ovtchinnikov}, {\em {J}acobi correction equation, line search, and
  conjugate gradients in {H}ermitian eigenvalue computation {I}: computing an
  extreme eigenvalue}, SIAM J. Numer. Anal., 46 (2008), pp.~2567--2592.

\bibitem{ovtchinnikov2008jacobiII}
{\sc E.~E. Ovtchinnikov}, {\em {J}acobi correction equation, line search, and
  conjugate gradients in {H}ermitian eigenvalue computation {II}: computing
  several extreme eigenvalues}, SIAM J. Numer. Anal., 46 (2008),
  pp.~2593--2619.

\bibitem{perdonG1986exterme}
{\sc A.~Perdon and G.~Gambolati}, {\em Exterme eigenvalues of large sparse
  matrices by {Rayleigh} quotient and modified conjugate gradients}, Comput.
  Methods Appl. Mech. Engrg., 56 (1986), pp.~251--264.

\bibitem{polakR1969note}
{\sc E.~Polak and G.~Ribi\'ere}, {\em Note sur la convergence de m\'ethodes de
  directions conjugu\'ees}, Revue Fran\c{c}aise d'Informatique et de Recherche
  Op\'erationnelle, 16 (1969), pp.~35--43.

\bibitem{polyak1964some}
{\sc B.~Polyak}, {\em Some methods of speeding up the convergence of iteration
  methods}, USSR Computational Mathematics and Mathematical Physics, 4 (1964),
  pp.~1--17.

\bibitem{powell1976some}
{\sc M.~Powell}, {\em Some convergence properties of the conjugate gradient
  method}, Math. Program., 11 (1976), pp.~42--49.

\bibitem{ritter1980rate}
{\sc K.~Ritter}, {\em On the rate of superlinear convergencee of a class of
  variable metric methods}, Numer. Math., 35 (1980), pp.~293--313.

\bibitem{samokish1958steepest}
{\sc B.~Samokish}, {\em The steepest descent method for an eigenvalue problem
  with semi-bounded operators}, Izv. Vyssh. Uchebn. Zaved. Mat., 5 (1958),
  pp.~105--114.
\newblock in Russian.

\bibitem{takahashi1965note}
{\sc I.~Takahashi}, {\em A note on the conjugate gradient method}, Inform.
  Process. Japan, 5 (1965), pp.~45--49.

\bibitem{wuXS2019trplk}
{\sc L.~Wu, F.~Xue, and A.~Stathopoulos}, {\em {TRPL+K}: thick-restart
  preconditioned {Lanczos+K} method for large symmetric eigenvalue problems},
  SIAM J. Sci. Comput., 41 (2019), pp.~A1013--1040.

\bibitem{zhangX2024chebyshev}
{\sc T.~Zhang and F.~Xue}, {\em A {Chebyshev} locally optimal block
  preconditioned conjugate gradient method for product and standard symmetric
  eigenvalue problems}, SIAM J. Matrix Anal. Appl., 45 (2024), pp.~2211--2242.

\end{thebibliography}

\appendix
\section{Exact line search in the block version}\label{sec:appendix}
\begin{theorem} \label{thm:linesearch}
	Let $X,r(X),V\in \mathbb{C}^{n\times n_b}$ be of full column rank, and $W\in \mathbb{C}^{n\times n_bm}$, which satisfy $V^{\HH}r(X)\ne0, W^{\HH}r(X)=0$.
	Suppose that $[X, V,W] $ has full column rank,
	and
	$( a_+, b_+)$ is a stationary point of the function $\tr\rho\left(X+ (I-P(X))[Va+Wb] \right)$.
	Write
	\[
		D =   (I-P(X))(V a_++Wb_+),\quad
		X_+ = X + D.
	\]
	Then, for the nontrivial case that $a_+$ is nonsingular, 
	\begin{subequations}
	\begin{align}
		\label{eq:thm:linesearch:r_opt-perp}
		&\quad
		r(X_+)\perp\range([X,V,W,D]),
		\\
		a_+ 
		&= -(\wtd r(X)^{\HH}V)^{-1}\opL_{X_+}(D)^{\HH}D, \label{eq:thm:linesearch:alpha}
		\\ \wtd \rho(X_+)-\wtd \rho(X) &= (X^{\HH}X)^{-1}\wtd r(X)^{\HH}Va_+ =-(X^{\HH}X)^{-1}\opL_{X_+}(D)^{\HH}D,
		\label{eq:thm:linesearch:rho}
		\\ \wtd r(X_+)-\wtd r(X) &=(I-P(X))\opL_{X_+}(D),  \label{eq:thm:linesearch:r}
		\\ \opL_{X_+;(I-P(X))W}(b_+)&= -W^{\HH}(I-P(X))\opL_{X_+}(Va_+).
		\label{eq:thm:linesearch:beta}
	\end{align}
	\end{subequations}
Here $\wtd\rho,\wtd r,\opL$ are defined as follows: for $X,Y\in \mathbb{C}^{n\times n_b},T\in \mathbb{C}^{n\times n_bm},Z\in \mathbb{C}^{n_bm\times n_b}$,
\begin{align*}
	\wtd \rho(X)&=(X^{\HH}X)^{-1}X^{\HH}AX=(X^{\HH}X)^{-1/2}\rho(X)(X^{\HH}X)^{1/2},\\
	\wtd r(X)&=AX-X\wtd \rho(X)=r(X)(X^{\HH}X)^{1/2},\\
	\opL_{X}(Y)&=AY-Y\wtd\rho(X),\\
	\opL_{X;T}(Z)&=T^{\HH}ATZ-T^{\HH}TZ\wtd\rho(X).
\end{align*}
\end{theorem}
If $a_+$ is singular, then for any vector $y$ with $a_+y=0$, $X_+y=Xy$, which means
that $Xy$ has to be an eigenvector.
When this case occurs in the iterations of the algorithm,
$Xy$ will be kept without any change, or equivalently, $Xy$ is ``locked''.
Therefore, we may consider the iteration in the subspace $\range(A)$, where the new ``$X$'' has less columns but the new ``$a_+$'' is nonsingular.
Hence this case can be easily treated in the algorithm and its analysis.

\begin{proof}[Proof of~\cref{thm:linesearch}]
	First taking differentials gives
	$\diff\,\tr\rho(X)=2\Re\inner{\diff X,r(X)(X^{\HH}X)^{-1/2}}$, where $\inner{X,Y}=\tr(X^{\HH}Y)$ is the standard inner product of the unitary space $\mathbb{C}^{n\times n_b}$.
	In fact, it is evident that $\tr\rho(X)$ is differentiable with respect to $X$, thus its differential can be obtained from perturbation (``$\doteq$'' represents omitting order-$2$ terms $O(\N{Y}^2)$):
	\begin{align*}
		\MoveEqLeft[1]\tr\rho(X+ Y)
		\\&\doteq\tr\left([(X+ Y)^{\HH}(X+ Y)]^{-1}(X+ Y)^{\HH}A(X+ Y)\right)
		\\&\doteq\tr\left((I+(X^{\HH}X)^{-1}[X^{\HH} Y+ Y^{\HH}X])^{-1}(X^{\HH}X)^{-1}(X^{\HH}AX+X^{\HH}A Y+ Y^{\HH}AX))\right)
		\\&\doteq\tr\left([I-(X^{\HH}X)^{-1}[X^{\HH} Y+ Y^{\HH}X]](\wtd\rho(X)+(X^{\HH}X)^{-1}[X^{\HH}A Y+ Y^{\HH}AX])\right)
		\\&\doteq\tr\left(\wtd\rho(X)-(X^{\HH}X)^{-1}[X^{\HH} Y+ Y^{\HH}X])\wtd\rho(X)+(X^{\HH}X)^{-1}[X^{\HH}A Y+ Y^{\HH}AX]\right)
		\\&\doteq\tr\left(\wtd\rho(X)+(X^{\HH}X)^{-1}[X^{\HH}A Y-X^{\HH} Y\wtd\rho(X)+ Y^{\HH}\wtd r(X)]\right)
		\\&\doteq\tr\rho(X)+2\Re\inner{ Y, r(X)(X^{\HH}X)^{-1/2}}
		.
	\end{align*}

	Since $( a_+, b_+)$ is a stationary point of the function $\tr\rho$,
		\begin{align*}
			0 &= \diff\,\tr\rho\left(X+ (I-P(X))(Va_++ Wb_+) \right)
			\\&= 2\Re\inner{\diff\,[X+ (I-P(X))(Va+ Wb) ] , r(X_+)(X_+^{\HH}X_+)^{-1/2}}
			\\&= 2\Re\inner{(I-P(X))(V\diff a+W\diff b), r(X_+)(X_+^{\HH}X_+)^{-1/2}}
			\\&= 2\Re\inner{\diff a, V^{\HH}(I-P(X))r(X_+)(X_+^{\HH}X_+)^{-1/2}}
			\\&\qquad+2\Re\inner{\diff b, W^{\HH}(I-P(X))r(X_+)(X_+^{\HH}X_+)^{-1/2}}
			,
		\end{align*}
	which implies 
		$V^{\HH}(I-P(X))r(X_+)=W^{\HH}(I-P(X))r(X_+)=0$.
	Then \cref{eq:thm:linesearch:r_opt-perp} holds.
	Also,
	\begin{align*}
		0 &=D^{\HH}(AX_+-X_+\wtd\rho(X_+))
		\\&=D^{\HH}(AX-X\wtd\rho(X)+X[\wtd\rho(X)-\wtd\rho(X_+)]+AD-D\wtd\rho(X_+))
		\\&=D^{\HH}\wtd r(X)+D^{\HH}\opL_{X_+}(D)
		\\&=(Va_++Wb_+)^{\HH}(I-P(X))\wtd r(X)+D^{\HH}\opL_{X_+}(D)
		\\&=a_+^{\HH}V^{\HH}\wtd r(X)+D^{\HH}\opL_{X_+}(D)
		,
	\end{align*}
	which gives $a_+=-(\wtd r(X)^{\HH}V)^{-1}\opL_{X_+}(D)^{\HH}D$, namely \cref{eq:thm:linesearch:alpha}.
	Since
		\begin{align*} 
			\wtd r(X_+)-\wtd r(X)
			&=AX_+-X_+\wtd \rho(X_+)-AX+X\wtd\rho(X)
			\\&= AD-D\wtd \rho(X_+)-X[\wtd\rho(X_+)-\wtd\rho(X)]
			\\&= \opL_{X_+}(D)-X[\wtd\rho(X_+)-\wtd\rho(X)]
			,
		\end{align*}
		\cref{eq:thm:linesearch:r} holds because
		\[
					\wtd r(X_+)-\wtd r(X)
					=(I-P(X))(\wtd r(X_+)-\wtd r(X))
					=(I-P(X))\opL_{X_+}(D).
		\]
	Hence we have $X[\wtd\rho(X_+)-\wtd\rho(X)]=P(X)\opL_{X_+}(D)$, or equivalently,
	\begin{align*} \label{eq:thm:linesearch:prf:diff-rho}
			\wtd\rho(X_+)-\wtd\rho(X)
			&= (X^{\HH}X)^{-1}X^{\HH}[AD-D\wtd\rho(X_+)]
			\\&=  (X^{\HH}X)^{-1}X^{\HH}A(I-P(X))(Va_++Wb_+)\nonumber
			\\&=  (X^{\HH}X)^{-1}[AX-X\wtd\rho(X)]^{\HH}(I-P(X))(Va_++Wb_+)\nonumber
			\\&=  (X^{\HH}X)^{-1}\wtd r(X)^{\HH}Va_+\nonumber
			.
	\end{align*}
	Together with \cref{eq:thm:linesearch:alpha}, we obtain \cref{eq:thm:linesearch:rho}. 
	Finally, 
	\begin{align*}
		0 &= W^{\HH}(I-P(X))(r(X_+)-r(X)) 
		\\&= W^{\HH}(I-P(X))\opL_{X_+}(D)
		\\&= W^{\HH}(I-P(X))[A(I-P(X))(Va_++Wb_+)-(I-P(X))(Va_++Wb_+)\wtd\rho(X_+)]
		\intertext{
	by $W^{\HH}(I-P(X))AX
	=W^{\HH}AX-W^{\HH}\wtd \rho(X)=W^{\HH}\wtd r(X)=0$,
	}
		&= W^{\HH}(I-P(X))[A(Va_++Wb_+)-(Va_++Wb_+)\wtd\rho(X_+)]
		\\&= W^{\HH}(I-P(X))[AVa_+-Va_+\wtd\rho(X_+)]
		+ W^{\HH}(I-P(X))[AWb_+-Wb_+\wtd\rho(X_+)]
		\\&= W^{\HH}(I-P(X))\opL_{X_+}(Va_+)
		+ \opL_{X_+;(I-P(X))W}(b_+)
		,
	\end{align*}
	which implies \cref{eq:thm:linesearch:beta}.
\end{proof}

\section{Proofs}\label{sec:proofs}
\begin{proof}[{Proof of~\cref{eq:beta}}]
	By \cref{cor:-cite-theorem-2-2-bennerl2022convergence-}, 
		\begin{equation}\label{eq:alpha-beta}
			\alpha= \frac{d^{\HH}F_+d}{-r^{\HH}r},\qquad
			\beta=-\alpha\frac{x_-^{\HH}\check{F}_+r}{x_-^{\HH}\check{F}_+x_-}.
		\end{equation}
	Also by \cref{cor:-cite-theorem-2-2-bennerl2022convergence-},
	$
	d ^{\HH}F_+d =\delta ,
	\check{F}_+d =r_+-r .
	$
	Therefore, $\alpha 
	=\frac{\delta }{-r ^{\HH}r }.$ Since
	\begin{align*}
		x_-^{\HH}\check{F}_+r &=x_-^{\HH}(I-P )F_+(I-P )r \\
		&=x_-^{\HH}F_+r -x_-^{\HH}x \cdot x ^{\HH}F_+r \\
		&=x_-^{\HH}(F_--\rho_+I+\rho_-I)r -\N{x_-}^2 x ^{\HH}(F -\rho_+I+\rho I)r \\
		&
		=\N{r }^2\cos^2\varphi_-
		,
		\\
		x_-^{\HH}\check{F}_+x_-&=x_-^{\HH}(I-P )F_+(I-P )x_-\\
		&=(x_--x \cos^2\varphi_-)^{\HH}F_+(x_--x \cos^2\varphi_-)\\
		&=x_-^{\HH}F_+x_--2\Re (x ^{\HH}F_+x_-)\cos^2\varphi_-+ x ^{\HH}F_+x \cos^4\varphi_-
		\\
		&=(\delta_-+\delta )\cos^2\varphi_--2\delta \cos^4\varphi_-+\delta \cos^4\varphi_-\\
		&=(\delta_-+\delta \sin^2\varphi_-)\cos^2\varphi_-.
	\end{align*}
	Then the result follows directly.
\end{proof}

\begin{proof}[{Proof of \cref{eq:comparison}}]
As preparation,
calculation gives
\begin{align} \label{eq:preparations}
	\inner{x ,d }_*&=(r +\varepsilon x )^{\HH}d =r ^{\HH}d =\alpha r ^{\HH}r =-\delta ,\\
\nonumber	\N{d }_*^2
	&=\inner{x_+,x_+}_*-\inner{x ,d }_*-\inner{d ,x }_*-\inner{x ,x }_*
	=\delta +\varepsilon_+\tan^2\varphi 
	,
	\\
\nonumber	\inner{r ,(I-P )x_-}_*
	&=\inner{r ,x_--x \cos^2\varphi_-}_*
	\\\nonumber&=r ^{\HH}\left(r_-+\varepsilon_-x_--(r +\varepsilon x )\cos^2\varphi_-\right)
	=-\N{r }^2\cos^2\varphi_-
	,
	\\\nonumber
	\N{(I-P )x_-}^2
	&=\N{x_--x \cos^2\varphi_-}^2
	=\cos^2\varphi_-\sin^2\varphi_-
	.
	\\\nonumber
	\N{(I-P )x_-}_*^2
	&=x_-^{\HH}\check F_+x_-+\varepsilon_+\N{(I-P )x_-}^2
	\\\nonumber&=(\delta_-+\delta \sin^2\varphi_-)\cos^2\varphi_-+\varepsilon_+\cos^2\varphi_-\sin^2\varphi_-
	\\\nonumber&=(\delta_-+\varepsilon \sin^2\varphi_-)\cos^2\varphi_-
	\\\nonumber
	\N{d }_*^2
	&=\N{\alpha r +\beta (I-P )x_-}_*^2
	\\\nonumber&=\N{\alpha r }_*^2-2\beta \alpha \N{r }^2\cos^2\varphi_-+\beta ^2(\delta_-+\varepsilon \sin^2\varphi_-)\cos^2\varphi_-
	\\\nonumber&=\N{\alpha r }_*^2+\beta \cos^2\varphi_-\left(2\delta +\beta (\delta_-+\varepsilon \sin^2\varphi_-)\right)
	\\\nonumber&=\N{\alpha r }_*^2-\frac{\delta ^2\cos^2\varphi_-({\delta_-+(\delta -\varepsilon_+)\sin^2\varphi_-})}{(\delta_-+\delta \sin^2\varphi_-)^2}
	,
\end{align}
	which implies
	\begin{align*}
		\N{\alpha r }_*^2
		&=\frac{\delta ^2\cos^2\varphi_-({\delta_-+(\delta -\varepsilon_+)\sin^2\varphi_-})}{(\delta_-+\delta \sin^2\varphi_-)^2}
	+\delta +\varepsilon_+\tan^2\varphi 
	\\&=\frac{\delta ^2\cos^2\varphi_-}{\delta_-+\delta \sin^2\varphi_-}+\delta -\alpha \delta \varepsilon_+
	,
	\end{align*}
	by $\tan^2\varphi =d ^{\HH}d 
		=\N{\beta (I-P )d_-}^2+\N{\alpha r }^2
		=\beta ^2\cos^2\varphi_-\sin^2\varphi_--\alpha \delta 
	$ and \cref{eq:beta}.
Thus,
\begin{align*}\nonumber
	\MoveEqLeft[1]\N{x +\alpha (A-\lambda_1I)x }_*^2
	\\\nonumber&=(x +\alpha (r +\varepsilon x ))^{\HH}(A-\lambda_1I)(x +\alpha (r +\varepsilon x ))\\\nonumber
	&=((1+\alpha \varepsilon )x +\alpha r )^{\HH}(F +\varepsilon I)((1+\alpha \varepsilon )x +\alpha r )\\\nonumber
	&=(2\alpha \varepsilon_++\alpha ^2\varepsilon ^2)\varepsilon +\varepsilon -2\delta +\N{\alpha r }_*^2\\\nonumber
	&=(2\alpha \varepsilon_++\alpha ^2\varepsilon ^2)\varepsilon +\varepsilon -2\delta +
	\frac{\delta ^2\cos^2\varphi_-}{\delta_-+\delta \sin^2\varphi_-}+\delta -\alpha \delta \varepsilon_+
	\\&=(\alpha \varepsilon_++\alpha ^2\varepsilon ^2)\varepsilon +\alpha \varepsilon_+^2
	+\varepsilon_+
	+\frac{\delta ^2}{\delta_-+(\delta_-+\delta )\tan^2\varphi_-}
	,
\end{align*}
which is exactly \cref{eq:comparison}.
\end{proof}

\begin{proof}[{Proof of \cref{eq:get-rid-of-phi}}]
	Since $\rho \le \rho_-\le \rho_0< \lambda_2$, we can write $\frac{\varepsilon_0}{\lambda_2-\lambda_1-\varepsilon_0}=\tau_0^2$.
	Then letting $x =\sum\limits_{i=1}^n\xi_{(i)}u_i$, $\varepsilon =\sum\limits_{i\ge 2}\abs{\xi_{(i)}}^2(\lambda_i-\lambda_1)\ge (\lambda_2-\lambda_1)\sin^2\angle(x ,u_1)$,
	$\N{r }^2=\sum\limits_{i\ge 2}\abs{\xi_{(i)}}^2(\lambda_i-\lambda_1)(\lambda_i-\rho )\in\left[\frac{\lambda_2-\lambda_1}{1+\tau_0^2}\varepsilon ,(\lambda_n-\lambda_1)\varepsilon \right]$.
	Also $\frac{\varepsilon }{\lambda_2-\lambda_1-\varepsilon }\ge \tan^2\angle(x ,u_1)$,
	$\varepsilon_0\le \frac{\tau_0^2(\lambda_2-\lambda_1)}{1+\tau_0^2}$.
	Thus $\varepsilon \ge \delta =\abs{\alpha }\N{r }^2\ge\frac{\lambda_2-\lambda_1}{1+\tau_0^2}\abs{\alpha }\varepsilon $,
	which implies $\abs{\alpha }\le \frac{1+\tau_0^2}{\lambda_2-\lambda_1}$.

	Note that 
	$\left({1+\frac{\delta_-+\delta }{\delta_-}\tan^2\varphi_-}\right)^{-1}\ge 1-\frac{\delta_-+\delta }{\delta_-}\sin^2\varphi_-$,
	and
	\begin{align}\nonumber
		\sin^2 \varphi_-
		&\le (\sin\angle(x ,u_1)+\sin\angle(x_-,u_1))^2
		\\&\le 2\left(\sin^2\angle(x ,u_1)+\sin^2\angle(x_-,u_1)\right)
		\le\frac{2(\varepsilon +\varepsilon_-)}{(\lambda_2-\lambda_1)^2}
		\le\frac{4\varepsilon_-}{(\lambda_2-\lambda_1)^2}
		.
		\label{eq:sin}
	\end{align}
	Under the fact $\delta \ge (1-C^{\SD})\varepsilon$,
	\begin{align*}
		\frac{\delta ^2}{\delta_-+(\delta_-+\delta )\tan^2\varphi_-}
		&\ge 
		\frac{\delta ^2}{\delta_-}\left(1-\frac{\delta_-+\delta }{\delta_-}\sin^2\varphi_-\right)
		\\&\ge 
		\frac{\delta ^2}{\delta_-}-\frac{\delta ^2}{(1-C^{\SD})\varepsilon_-}\frac{\varepsilon_-}{(1-C^{\SD})\varepsilon_-}\frac{4\varepsilon_-}{(\lambda_2-\lambda_1)^2}
		\\&\ge 
		\frac{\delta ^2}{\delta_-}-\frac{4\delta ^2}{(1-C^{\SD})^2(\lambda_2-\lambda_1)^2}
		,
	\end{align*}
	which guarantees \cref{eq:get-rid-of-phi}.
\end{proof}

\begin{proof}[{Proof of \cref{thm:convergence-rate:1m1}}]
By \cref{cor:-cite-theorem-2-2-bennerl2022convergence-}, 
we still have \cref{eq:alpha-beta} and also
$\begin{bmatrix}
	\beta_k\\\alpha_kc_k
\end{bmatrix}=-(W_k^{\HH}\check{F}_{k+1}W_k)^{-1}W_k^{\HH}\check{F}_{k+1}r_k.$
All the equalities in the proof of \cref{eq:beta} remain valid.

Similarly we simplify the subscripts related to $k$.
	By \cite{golubY2002inverse}, The eSD method $\LOCG(1,m_e,0)$  enjoys
	\begin{equation*}\label{eq:rate:eSD}
		\varepsilon_+\le C^{\SD}\varepsilon ,
		\quad C^{\SD}= 
		{\cheb_{m_e}(\Delta^{-1})}^{-2}+O(\sqrt{\varepsilon })<1,
	\end{equation*}
	where $\cheb_{m_e}(t)$ is the $m_e$-th degree Chebyshev polynomial of the first kind.
Similarly, we compare the convergence behaviors of $\LOCG(1,m_e,1)$ and eSD method $\LOCG(1,m_e,0)$.
We aim to relate $\N{x_+}_*^2=\varepsilon_+$ and $\N{x +\alpha q }_*^2$, where $q =r +(I-P )(I-Q )W c $. 
Note that \cref{eq:preparations} still holds, except the last one becomes
\begin{equation*} 
\begin{aligned}[t]
	\N{d }_*^2
	&=\N{\alpha q +\beta (I-P )x_-}_*^2
	\\&=\N{\alpha q }_*^2+2\beta \alpha \Re\inner{q ,(I-P )x_-}_*+\beta ^2(\delta_-+\varepsilon \sin^2\varphi_-)\cos^2\varphi_-
	.
\end{aligned}
\end{equation*}
Since 
\begin{align*}
	\inner{(I-P )x_-,q -r }_*
	&=\inner{(I-P )x_-,(I-P )(I-Q )W c }_*
	\\&=(x_--x \cos^2\varphi_-)^{\HH}(A-\lambda_1I)(I-P -Q )W c 
	\\&=([r_-+\varepsilon_-x_-]-[r +\varepsilon x ]\cos^2\varphi_-)^{\HH}(I-P -Q )W c 
	\\&=-[r_-+\varepsilon_-x_-]^{\HH}P W c 
	\\&=-[x_-^{\HH}F_-x +\varepsilon_-x_-^{\HH}x ]x ^{\HH}W c 
	\\&=-\varepsilon \cos^2\varphi_-\N{r }^2\tau ,
\end{align*}
where $\tau =\frac{x ^{\HH}W c }{r ^{\HH}r }$,
we have
\[
	\N{d }_*^2=
	\N{\alpha q }_*^2-\frac{\delta ^2\cos^2\varphi_-\varepsilon \tau }{\delta_-+\delta\sin^2\varphi_-}-\frac{\delta ^2\cos^2\varphi_-({\delta_-+(\delta -\varepsilon_+)\sin^2\varphi_-})}{(\delta_-+\delta \sin^2\varphi_-)^2},
\]
and then
\[
		\N{\alpha q }_*^2
	=
	\frac{(1+\varepsilon \tau )\delta ^2\cos^2\varphi_-}{\delta_-+\delta \sin^2\varphi_-}+\delta -\alpha \delta \varepsilon_+
	.
\]
Substituting all the $F $ in $q $ with $A-\lambda_1I$ to get $\wtd q $, then $\wtd q =q +\varepsilon x +O(\varepsilon ^{3/2})$, and comparing with \cref{eq:comparison},
\begin{align*}
		\N{x +\alpha \wtd q }_*^2
		&=\N{x +\alpha \left(q +\varepsilon x +O(\varepsilon ^{3/2})\right)}_*^2
		\\&=
		\varepsilon_+
		+\frac{\delta ^2(1+\varepsilon \tau )}{\delta_-+(\delta_-+\delta )\tan^2\varphi_-}
		+O(\varepsilon ^2)
		.
	\end{align*}
Expressing $W ,c $ in $r , F $, it is not difficult to show $\tau =O(\N{r }^2)=O(\delta )$.
Using the same argument in \cref{sec:a-simple-case},
		$\N{x +\alpha \wtd q }_*^2
		= \varepsilon_+ +\frac{\delta ^2}{\delta_-}
		+O(\varepsilon ^2).$
Going on, we arrive at \cref{thm:convergence-rate:1m1}.
\end{proof}

\begin{proof}[{Proof of \cref{thm:convergence-rate:11m}}]
By \cref{cor:-cite-theorem-2-2-bennerl2022convergence-}, 
we still have $\alpha_k=\frac{\delta_k}{-r_k^{\HH}r_k}$ and also
$b_k=-\alpha_k(W_k^{\HH}\check{F}_{k+1}W_k)^{-1}W_k^{\HH}\check{F}_{k+1}r_k.$

Again we simplify the subscripts related to $k$, but this time $\cdot_{k-i}$ is replaced by $\cdot_{-i}$.
Similarly as the discussion in \cref{sec:a-simple-case}, we calculate a relation between $\N{x_+}_*^2$ and $\N{x +\alpha (A-\lambda_1I)x }_*^2$.
Write 
\[
	e =\frac{W ^{\HH}F r }{r ^{\HH}r }=\frac{1}{r ^{\HH}r }\begin{bmatrix}
	r ^{\HH}r_{{-1}},\dots,r ^{\HH}r_{{-m_h}}
\end{bmatrix}^{\HH},f =W ^{\HH}x =\begin{bmatrix}
x ^{\HH}x_{{-1}},\dots,x ^{\HH}x_{{-m_h}}
\end{bmatrix}^{\HH}.
\]
As an analogue of \cref{eq:preparations}, we have
\begin{align}\label{eq:quantities:1}
	b &=\delta (W ^{\HH}\check{F}_+W )^{-1}(e -f ),\\\nonumber
b ^{\HH}W ^{\HH}\check F_+W b 
&=\delta (e -f )^{\HH}b \in \mathbb{R},
	\\\nonumber
	\inner{r ,(I-P )W b }_*
	&=\N{r }^2(e -f )^{\HH}b  \in \mathbb{R},
	\\\nonumber
	\N{(I-P )W b }_*^2
	&=\delta (e -f )^{\HH}b +\varepsilon_+\N{W b }^2-\varepsilon_+\abs{f ^{\HH}b }^2
	,
	\\\nonumber
	\tan^2\varphi 
	&=d ^{\HH}d 
	=\N{W b }^2-\abs{f ^{\HH}b }^2-\alpha \delta ,
		\\\nonumber
	\delta +\varepsilon_+\tan^2\varphi 
	&=\N{d }_*^2
	=\N{\alpha r +(I-P )W b }_*^2
	\\\nonumber&=\N{\alpha r }_*^2-\delta (e -f )^{\HH}b +\varepsilon_+(\N{W b }^2-\abs{f ^{\HH}b }^2)
	,
\end{align}
	which implies 
	\begin{equation}\label{eq:quantities:2}
		\N{\alpha r }_*^2
		=
	\delta -\varepsilon_+\alpha \delta 
	+\delta (e -f )^{\HH}b 
	.
	\end{equation}
	Thus, as an analogue of \cref{eq:comparison},
	\begin{align}\nonumber
		\MoveEqLeft \N{x +\alpha (A-\lambda_1I)x }_*^2
		\\\nonumber&=(2\alpha \varepsilon_++\alpha ^2\varepsilon ^2)\varepsilon +\varepsilon -2\delta +\N{\alpha r }_*^2\\\nonumber
		&=(\alpha \varepsilon_++\alpha ^2\varepsilon ^2)\varepsilon +\alpha \varepsilon_+^2
		+\varepsilon_+
	+\delta (e -f )^{\HH}b 
		\\&=O(\varepsilon ^2) +\varepsilon_+ +\delta ^2(e -f )^{\HH}(W ^{\HH}\check{F}_+W )^{-1}(e -f )
		.
		\label{eq:quantities:3}
	\end{align}
	The rest is to estimate $(W ^{\HH}\check{F}_+W )^{-1}$ and $e -f $.
	Since
	\begin{equation}\label{eq:quantities:4}
	W ^{\HH}\check{F}_+W 
	=(W ^{\HH}-f x ^{\HH})(F_{k}+\delta I)(W -x f ^{\HH})
	=W ^{\HH}F_+W -\delta f f ^{\HH}
	,
	\end{equation}
whose $(i,j)$th entry for $i\le j$ is
\begin{align*}
			\MoveEqLeft x_{{-i}}^{\HH} F_+x_{{-j}}-\delta x_{{-i}}^{\HH}x x ^{\HH}x_{{-j}}
			\\&=x_{{-i}}^{\HH}(F_{{-i}}+(\varepsilon_{{-i}}-\varepsilon_+)I)x_{{-j}}-\delta x_{{-i}}^{\HH}x x ^{\HH}x_{{-j}}
			\\&=(\varepsilon_{{-i}}-\varepsilon_+)x_{{-i}}^{\HH}x_{{-j}}-\delta x_{{-i}}^{\HH}x x ^{\HH}x_{{-j}}
			\\&=\left[(\varepsilon_{{-i}}-\varepsilon_+)\cos\varphi_{k-i,k-j}-\delta \right]\cos\varphi_{k,k-i}\cos\varphi_{k,k-j}(1+O(\varepsilon_{-i})+O(\varepsilon_{-j}))
			,
\end{align*}
	where $\varphi_{k-i,k-j}=\angle(x_{{-i}},x_{{-j}})$.
	Since $1-\cos\varphi_{k-i,k-j}=2\sin^2\frac{\varphi_{k-i,k-j}}{2}=O(\varepsilon_{{-j}})$ as is shown in \cref{eq:sin},
	\begin{align*}
		W ^{\HH}\check{F}_+W 
		&=[\varepsilon_{{-i}}-\varepsilon_+-\delta +O(\varepsilon_{{-m_h}})^2]
		\\&=\begin{bmatrix}
			\varepsilon_{{-1}}-\varepsilon  & \varepsilon_{{-1}}-\varepsilon  &\dots & \varepsilon_{{-1}}-\varepsilon \\
			\varepsilon_{{-1}}-\varepsilon  & \varepsilon_{{-2}}-\varepsilon  &\dots & \varepsilon_{{-2}}-\varepsilon \\
			\vdots & \vdots &\ddots &  \\
			\varepsilon_{{-1}}-\varepsilon  & \varepsilon_{{-2}}-\varepsilon  & & \varepsilon_{{-m_h}}-\varepsilon \\
		\end{bmatrix}
		+O(\varepsilon_{{-m_h}}^2)
		\\&=J^{-1}\diag(\delta_{{-1}},\delta_{{-2}},\dots,\delta_{{-m_h}})J^{-\HH}
		+O(\varepsilon_{{-m_h}}^2)
	\end{align*}
	where $J=\begin{bmatrix}
		1 & &&\\
		-1 & 1&&\\
		&\ddots&\ddots&\\
		&&-1&1\\
	\end{bmatrix}$. 
	Hence 
	\[
		(W ^{\HH}\check{F}_+W )^{-1}
		=J^{\HH}\diag(\delta_{{-1}}^{-1},\delta_{{-2}}^{-1},\dots,\delta_{{-m_h}}^{-1})J\left(I+ O(\delta_{{-1}}^{-1}\varepsilon_{{-m_h}}^2)\right).
	\]
	Since 
	$\abs{d_{{-i}}^{\HH}x }
		=\abs{d_{{-i}}^{\HH}(x_{{-i}}+d_{{-i}}+\dots+d_{{-1}})}
		\le \N{x_{{-i}}}^2\sin^2\varphi_{{-i}}+\sum_{j=1}^i\N{d_{{-i}}}\N{d_{{-j}}}
		=O(\varepsilon_{{-i}})
		$, and then
	\begin{align*}
		J(e -f )
		&=
		\frac{1}{r ^{\HH}r }\begin{bmatrix}
			(r_{{-1}}-r )^{\HH}r +r ^{\HH}r \\ (r_{{-2}}-r_{{-1}})^{\HH}r \\\vdots\\(r_{{-m_h}}-r_{{-m_h+1}})^{\HH}r 
\end{bmatrix}-\begin{bmatrix}
1-d_{{-1}}^{\HH}x \\-d_{{-2}}^{\HH}x \\\vdots\\-d_{{-m_h}}^{\HH}x 
\end{bmatrix}
\\&=\begin{bmatrix}
	\gamma_1\\ \gamma_2\\\vdots\\\gamma_{m_h}
\end{bmatrix}
+\begin{bmatrix}
	O(\varepsilon_{{-1}})\\O(\varepsilon_{{-2}})\\\vdots\\O(\varepsilon_{{-m_h}})
\end{bmatrix},
\quad \text{for}\; \gamma_j=\frac{(r_{{-j}}-r_{{-j+1}})^{\HH}r }{r ^{\HH}r }.
	\end{align*}
	Note that $\gamma_1=-1$ for $r\perp r_{-1}$.
	To sum up, as an analogue of \cref{eq:basic-relation}, letting $\eta =\frac{\varepsilon_+}{\varepsilon }$ and 
		 $\N{x +\alpha (A-\lambda_1I)x }_*^2= C\N{x }_*^2=C\varepsilon_{k}$, we have
		 \[
			 \varepsilon_++\gamma_1^2\frac{\delta ^2}{\delta_{{-1}}}+\gamma_2^2\frac{\delta ^2}{\delta_{{-2}}}+\cdots+\gamma_{m_h}^2\frac{\delta ^2}{\delta_{{-m_h}}}=C\varepsilon +O(\varepsilon ^2),
		 \]
		 and
		 \begin{equation}\label{eq:optim:m}
			 \begin{multlined}[b]
			 \eta +\left(\gamma_1^2\frac{\eta_{{-1}}}{1-\eta_{{-1}}}+\gamma^2_2\frac{\eta_{{-1}}\eta_{{-2}}}{1-\eta_{{-2}}}+\cdots+\gamma^2_{m_h}\frac{\eta_{{-1}}\cdots \eta_{{-m_h}}}{1-\eta_{{-m_h}}}\right)(1-\eta )^2
			 \\=C+O(\varepsilon )=:\wtd C.
			 \end{multlined}
		 \end{equation}
		 Solving the $m_h$-variable optimization problem 
		 \[
		 	\max\; \eta \dots\eta_{{-m_h}}, \qquad\text{s.t. \cref{eq:optim:m} holds}
		 \]
		 produces 
		 \[
		 	\eta \dots\eta_{{-m_h}}\le \chi(\sigma,\wtd C)
			, \quad \sigma=(\abs{\gamma_1}+\abs{\gamma_2}+\dots+\abs{\gamma_{m_h}})^2.
		 \]
		 Using the same argument in \cref{sec:a-simple-case}, we have \cref{thm:convergence-rate:11m}.
\end{proof}

\begin{proof}[{Proof of \cref{thm:convergence-rate:m11}}]
Again we simplify the subscripts related to $k$. Moreover, we also drop the subscript $\cdot_{(1)}$ when no confusion occurs and replace the subscript $\cdot_{(2:n_b)}$ with $\cdot_{(:)}$.
The process is the same as that in \cref{sec:more-historical-terms} except the definition of $W$.
Hence \cref{eq:quantities:1,eq:quantities:2,eq:quantities:3,eq:quantities:4} still hold except that the meaning of $e,f$ has been changed as follows:
\[
	e^{\HH}=\frac{ r^{\HH}FW}{r^{\HH}r}
	=\left[ 0 , \frac{r^{\HH}R_{(:)}}{r^{\HH}r} , \frac{r^{\HH}FR_{\perp}}{r^{\HH}r} \right]
	,
	f^{\HH}=x^{\HH}W
	=\left[ x^{\HH}X_- , 0 , 0\right]
	.
\]
Then,
\begin{align*}
	W^{\HH}\check F_+W
	&=W^{\HH}F_+W-\delta_{(1)}ff^{\HH}
	\\
	&=\begin{bmatrix}
		E+\delta_--\delta_{(1)}gg^{\HH} &  \Gamma_{(:)}^{\HH}E_{(:)}                         & 0\\
		E_{(:)}\Gamma_{(:)}             & E_{(:)}           & R_{\perp}^{\HH}R_{\perp}\\
		 0                                                        & R_{\perp}^{\HH}R_{\perp} & R_{\perp}^{\HH}F_+R_{\perp}\\
	\end{bmatrix}
	=\begin{bmatrix}
		\Omega_{11} & \Omega_{21}^{\HH}\\ \Omega_{21} & \Omega_{22}
	\end{bmatrix}
	,
\end{align*}
where $\Omega_{11}\in \mathbb{C}^{n_b\times n_b}$.
Since $(\rho_+,x_+)$ is the minimal Ritz pair of $A$ on the subspace $\range([x,r,W])$, 
supposing $\rho_{0,(1)}<\lambda_2$,
$[x,r,W]^{\HH}F_+[x,r,W]\succeq 0$ and has only one zero eigenvalue,
which results in $W^{\HH}\check F_+W\succ 0$.
Thus, $\Omega_{11},\Omega_{22}$ and their Schur complements are all positive definite.
By the identities
\begin{align*}
\begin{bmatrix}
	\Omega_{11} & \Omega_{21}^{\HH}\\ \Omega_{21} & \Omega_{22}
	\end{bmatrix}^{-1}
	&=(*)^{\HH}\begin{bmatrix}
\Omega_{11} & \\ & \Omega_{22}-\Omega_{21}\Omega_{11}^{-1}\Omega_{21}^{\HH}
\end{bmatrix}^{-1}\begin{bmatrix}
I &\\ -\Omega_{21}\Omega_{11}^{-1}& I
\end{bmatrix},
\end{align*}
in which $*$ denotes the terms determined by the Hermitian property,
\begin{align}
	\varsigma&=(e-f)^{\HH}(W^{\HH}\check F_+W)^{-1}(e-f)
	-
	g^{\HH}\Omega_{11}^{-1}g
	\nonumber\\ &=(*)^{\HH}
	\begin{bmatrix}
		E_{(:)}-E_{(:)}\Gamma_{(:)}\Omega_{11}^{-1}\Gamma_{(:)}^{\HH}E_{(:)}            & R_{\perp}^{\HH}R_{\perp} \\
		R_{\perp}^{\HH}R_{\perp} & R_{\perp}^{\HH}F_+R_{\perp}
	\end{bmatrix}^{-1}
	\begin{bmatrix}
E_{(:)}\Gamma_{(:)}\Omega_{11}^{-1}g
		+	\frac{R_{(:)}^{\HH}r}{r^{\HH}r} \\
		\frac{R_{\perp}^{\HH}Fr}{r^{\HH}r}
	\end{bmatrix}
	\nonumber\\ &=:(*)^{\HH}
	\begin{bmatrix}
		\Psi_{11} & \Psi_{21} \\
		\Psi_{21} & \Psi_{22}
	\end{bmatrix}^{-1}
	\begin{bmatrix}
		\Phi_0+\wtd \Phi_1\\
		\Phi_2\\
	\end{bmatrix}
	,\qquad \text{where}\;\Phi_0=E_{(:)}\Gamma_{(:)}\Omega_{11}^{-1}g
	\nonumber\\&=(*)^{\HH}[\Psi_{11}-\Psi_{21}\Psi_{22}^{-1}\Psi_{21}]^{-1}[\Phi_0+\Phi_1]+\Phi_2^{\HH}\Psi_{22}^{-1}\Phi_2
	,
	\label{eq:varsigma}
\end{align}
where $\Phi_1=\wtd\Phi_1-\Psi_{21}\Psi_{22}^{-1}\Phi_2$.
Here $\Psi_{11},\Psi_{22},\Psi_{11}-\Psi_{21}\Psi_{22}^{-1}\Psi_{21}$ are all positive definite.

The rest is to find out its dominating term.
By the discussion in \cref{rk:h:new-order}, the previous estimate of $\Gamma$ gives
\begin{align*}
	\Omega_{11}&=E+\delta_--\delta_{(1)}gg^{\HH}
	=\underbrace{\begin{bmatrix}
		\delta_{-(1)}&\\&E_{(:)}+\delta_{-(:)}
\end{bmatrix}}_{\Omega_\theta}+\delta_{(1)}(e_1e_1^{\HH}-gg^{\HH}),
\\[-2.7ex]
\Omega_{11}^{-1}
&=\Omega_\theta^{-1}-\delta_{(1)}\Omega_\theta^{-1}[e_1,g]{\overbrace{\left(\begin{bsmallmatrix}
1 & \\&-1
\end{bsmallmatrix}+\delta_{(1)}[e_1,g]^{\HH}\Omega_\theta^{-1}[e_1,g]\right)}^{\Omega_s}{}^{-1}}[e_1,g]^{\HH}\Omega_\theta^{-1}
\\ &=\Omega_\theta^{-1}-\frac{\delta_{(1)}}{\det(\Omega_s)}\bigg((\delta_{(1)}\omega_{gg}-1)\Omega_\theta^{-1}e_1e_1^{\HH}\Omega_\theta^{-1}-\delta_{(1)}\ol\omega_{ge}\Omega_\theta^{-1}e_1g^{\HH}\Omega_\theta^{-1}
\\&\qquad\qquad\qquad\qquad\qquad-\delta_{(1)}\omega_{ge}\Omega_\theta^{-1}ge_1^{\HH}\Omega_\theta^{-1}+(1+\delta_{(1)}\omega_{ee})\Omega_\theta^{-1}gg^{\HH}\Omega_\theta^{-1}\bigg)
,
\end{align*}
where 
\begin{align*}
	\omega_{ee} &=e_1^{\HH}\Omega_\theta^{-1}e_1=\delta_{-(1)}^{-1},\\
	\omega_{gg}&=g^{\HH}\Omega_\theta^{-1}g=\abs{ g_{(1)}}^2\delta_{-(1)}^{-1}+g_{(:)}^{\HH}(E_{(:)}+\delta_{-(:)})^{-1}g_{(:)}=\delta_{-(1)}^{-1}(1+O(\varepsilon_{(1)})),\\
	\omega_{ge}&=g^{\HH}\Omega_\theta^{-1}e_1= \ol g_{(1)}\delta_{-(1)}^{-1}=\delta_{-(1)}^{-1}(1+O(\varepsilon_{(1)})),\\
	\det(\Omega_s) &= (1+\delta_{(1)}\omega_{ee})(\delta_{(1)}\omega_{gg}-1)-\delta_{(1)}^2\abs{\omega_{ge}}^2
	\\ &= -1-\delta_{(1)}(\omega_{ee}-\omega_{gg})+\delta_{(1)}^2(\omega_{gg}\omega_{ee}-\abs{\omega_{ge}}^2)
		=-1+O(\N{\varepsilon_-}).
\end{align*}
Writing $\what g=[E_{(:)}+\delta_{-(:)}]^{-1/2}g_{(:)},\what H=H_{(:)}[E_{(:)}+\delta_{-(:)}]^{-1/2}$,
we have
$\omega_{gg}\omega_{ee}-\abs{\omega_{ge}}^2=\omega_{ee}\what g^{\HH}\what g$, and 
\begin{align*}
\Omega_{11}^{-1}g
&=
\Omega_\theta^{-1}g-\frac{\delta_{(1)}}{\det(\Omega_s)}\bigg((\delta_{(1)}\omega_{gg}-1)\ol\omega_{ge}\Omega_\theta^{-1}e_1-\delta_{(1)}\ol\omega_{ge}\omega_{gg}\Omega_\theta^{-1}e_1
\\&\qquad\qquad\qquad\qquad\qquad-\delta_{(1)}\abs{\omega_{ge}}^2\Omega_\theta^{-1}g+(1+\delta_{(1)}\omega_{ee})\omega_{gg}\Omega_\theta^{-1}g\bigg)
\\&=
\frac{(1+\delta_{(1)}\omega_{ee})\Omega_\theta^{-1}g-\delta_{(1)}\ol\omega_{ge}\Omega_\theta^{-1}e_1}{-\det(\Omega_s)}
,
\\
g^{\HH}\Omega_{11}^{-1}g
&=
\frac{\omega_{gg}+\delta_{(1)}\omega_{ee}\what g^{\HH}\what g}{-\det(\Omega_s)}
=
\delta_{-(1)}^{-1}(1+O(\N{\varepsilon_-})).
\end{align*}
Also,
\begin{align*}
	\Gamma_{(:)}\Omega_\theta^{-1}e_1&=\omega_{ee}h ,\qquad
	\Gamma_{(:)}\Omega_\theta^{-1}g=\ol\omega_{ge}h + \what H\what g,\qquad
	\Gamma_{(:)}\Omega_\theta^{-1}\Gamma_{(:)}^{\HH}=\delta_{-(1)}^{-1}h h ^{\HH}+ \what H\what H^{\HH},
	\\
\Gamma_{(:)}\Omega_{11}^{-1}g
&=
\Gamma_{(:)}\frac{(1+\delta_{(1)}\omega_{ee})\Omega_\theta^{-1}g-\delta_{(1)}\ol\omega_{ge}\Omega_\theta^{-1}e_1}{-\det(\Omega_s)}
\\ &
=\frac{\ol\omega_{ge}h +(1+\delta_{(1)}\omega_{ee})\what H\what g}{-\det(\Omega_s)}
=\delta_{-(1)}^{-1}h(1+O(\N{\varepsilon_-^{(1+\theta)/2}}))
,
\\
\Gamma_{(:)}\Omega_{11}^{-1}\Gamma_{(:)}^{\HH}
&=\omega_{ee}hh^{\HH}+\what H\what H^{\HH}
-\frac{\delta_{(1)}}{\det(\Omega_s)}\bigg((\delta_{(1)}\omega_{gg}-1)\omega_{ee}^2h h ^{\HH}
	-{\delta_{(1)}\ol\omega_{ge}}\omega_{ee}h(\ol\omega_{ge}h+\what H\what g)^{\HH} 
	\\ &\qquad\qquad \qquad
	-{\delta_{(1)}\omega_{ge}}\omega_{ee}(\ol\omega_{ge}h+\what H\what g)h^{\HH}
+(1+\delta_{(1)}\omega_{ee})(\ol\omega_{ge}h+\what H\what g)(\ol\omega_{ge}h+\what H\what g)^{\HH}
\bigg)
\\&=
\what H\what H^{\HH}-\frac{\omega_{ee}h h ^{\HH}
}{\det(\Omega_s)}+\frac{\delta_{(1)}}{\det(\Omega_s)}\bigg[\what g^{\HH}\what g\omega_{ee}h h ^{\HH}
	-\ol\omega_{ge}h\what g^{\HH}\what H^{\HH}
	-\omega_{ge}\what H\what gh^{\HH}
	 \\&\qquad\qquad \qquad \qquad\qquad \qquad\qquad \qquad
-(1+\delta_{(1)}\omega_{ee})\what H\what g\what g^{\HH}\what H^{\HH}\bigg]
\\&=(1+O(\N{\varepsilon_-^{(1+\theta)/2}}))\left(\delta_{-(1)}^{-1}h h ^{\HH}+H_{(:)}[E_{(:)}+\delta_{-(:)}]^{-1}H_{(:)}^{\HH}\right)
,
\\
\Phi_0^{\HH}\Psi_{11}^{-1}\Phi_0
&=(*)^{\HH}\left(E_{(:)}^{-1}-\Gamma_{(:)}\Omega_{11}^{-1}\Gamma_{(:)}^{\HH}\right)^{-1}\Gamma_{(:)}\Omega_{11}^{-1}g
\\&=(*)^{\HH}\Big(E_{(:)}^{-1}-\frac{hh^{\HH}}{\delta_{-(1)}}-H_{(:)}[E_{(:)}+\delta_{-(:)}]^{-1}H_{(:)}^{\HH}\Big)^{-1}\frac{h}{\delta_{-(1)}}(1+O(\N{\varepsilon_-^{(1+\theta)/2}}))
\\&=\frac{1+O(\N{\varepsilon_-^{(1+\theta)/2}})}{\delta_{-(1)}^2}h^{\HH}\Big(M-\frac{hh^{\HH}}{\delta_{-(1)}}\Big)^{-1}h
\\&=\frac{1+O(\N{\varepsilon_-^{(1+\theta)/2}})}{\delta_{-(1)}}\frac{\tau^2}{1-\tau^2}
,
\end{align*}
for $ \tau^2=\frac{h^{\HH}M^{-1}h}{\delta_{-(1)}}$.
However $\Psi_{11}\succ 0$, so $0<\tau^2<1$ and $M\succ 0$.
	Since 
	\[
		E_{(:)}^{-1}-H_{(:)}E_{(:)}^{-1}H_{(:)}^{\HH}=O(\N{\varepsilon_-^{(1+\theta)/2}}),\quad
	E_{(:)}^{-1}\delta_{-(:)}[E_{(:)}+\delta_{-(:)}]^{-1}=O(\N{\varepsilon_-}),
	\]
	we know 
	\[
		1>\tau^2=\frac{h^{\HH}M^{-1}h}{\delta_{-(1)}}\ge
		\frac{\N{h}^2}{\delta_{-(1)}\N{M}}
		\quad \Rightarrow \quad 
h=O(\delta_{-(1)}^{1/2}\N{\varepsilon_-^{(1+\theta)/4}}).
	\]
	This is \cref{eq:h:new-order}.

Return to \cref{eq:varsigma},
\begin{align*}
	\varsigma
	&=(1+O(\N{\varepsilon_-^{(1+\theta)/2}})[\delta_{-(1)}^{-1}h+\what\Phi_1]^{\HH}[M-\delta_{-(1)}^{-1}hh^{\HH}-\what\Psi_{21}\what\Psi_{22}\what\Psi_{21}^{\HH}]^{-1}[\delta_{-(1)}^{-1}h+\what\Phi_1]+\what\Phi_2^{\HH}\what\Psi_{22}\what\Phi_2
	,
\end{align*}
in which $\what\Psi_{22},\what\Psi_{21},\what\Phi_1,\what\Phi_2$ is as in \cref{eq:gamma:terms}.
Then we have
\begin{align*}
	\varsigma
	&=
	(1+O(\N{\varepsilon^{(1+\theta)/2}}))
	\frac{\gamma^2}{\delta_{-(1)}}
	,
\end{align*}
where $\gamma^2$ is as in \cref{eq:gamma:m11} or \cref{eq:gamma:m11:shrink}.
Therefore $(e-f)^{\HH}(W^{\HH}\check F_+W)^{-1}(e-f)=
g^{\HH}\Omega_{11}^{-1}g
+\varsigma
=
\frac{1+O(\N{\varepsilon_-^{(1+\theta)/2}})}{\delta_{-(1)}}(1+\gamma^2)
$,
and  finally we have
\[
	\varepsilon_{+(1)}+\frac{\delta_{(1)}^2}{\delta_{-(1)}}(1+\gamma^2)=C\varepsilon_{(1)}+O(\N{\varepsilon_-^{(3+\theta)/2}}),
\]
and 
\[
	\eta+\frac{\eta_-(1-\eta)^2}{1-\eta_-}(1+\gamma^2)
	=C+O(\N{\varepsilon_-^{(1+\theta)/2}})=:\wtd C.
\]
Going on, we have
		 \[
			 \eta\eta_-\le \chi(\sigma,\wtd C)
			 , \quad \sigma=1+\gamma^2.
		 \]
		 Using the same argument in \cref{sec:a-simple-case}, we have \cref{thm:convergence-rate:m11}.
\end{proof}
\end{document}